\documentclass[11pt,leqno,english]{article}

\usepackage{ShVa718_719}

\providecommand{\SimpleCite}[2][]{}

\newcommand{\FUnion} {\uplus}
\newcommand{\Tuple}  {\SimpleSeq}

\renewcommand{\Pair}[3][]{%
	\SimpleSeq[#1]{#2,#3}%
}

\newcommand{\Interval}[3][]{\MathX{%
	\IfEmpty{#1}%
		{\Braces{#2,\dots,#3}}%
	{%
		\IfEqual{#1}{+}%
			{\Braces{#2+1,\dots,#3}}%
    		{
			\IfEqual{#1}{-}%
				{\Braces{#2,\dots,#3-1}}%
				{\Braces{#2+1,\dots,#3-1}}%
		}%
	}%
}}

\providecommand{\Lan}[1]{\MathX{%
	{L}_{\infty#1}%
}}
\newcommand{\MLan}[2]{\MathX{%
	{M}_{\infty#1;#2}%
}}

\providecommand{\LEquiv}[1]{%
	\mathrel{{\equiv}_{\infty#1}}%
}
\newcommand{\MEquiv}[2]{%
	\mathrel{{\equiv}_{\infty#1;#2}}%
}

\newcommand{\Str}{\mathcal}

\newcommand{\EF}[5][]{\MathX{%
	\Symb{EF}_{#2;#3}\Par[#1]{#4,#5}%
}}

\newcommand{\PlayerOne}{\MathX{\forall}}
\newcommand{\PlayerTwo}{\MathX{\exists}}

\providecommand{\EhrFra}{Ehrenfeucht-Fra\"\i ss\'e}

\newcommand{\Param}{\MathX{P}}

\newcommand{\Equiv}[1][\phi,\Param]{%
	\mathrel{{\sim}_{\scriptscriptstyle #1}}%
}

\newcommand{\PSup}[2][]{\MathX{%
	{\sup}^{+}\Par[#1]{#2}%
}}

\newcommand{\Image}[3][]{\MathX{%
	#2\Brackets[#1]{#3}%
}}
\newcommand{\Proj}[3][]{\MathX{%
	#2\Brackets[#1]{#3}%
}}

\newcommand{\Closure}[3][]{\MathX{%
	\IfEmpty{#2}%
		{\MathSymb{Cl}\Par[#1]{#3}}
		{\MathSymb{Cl}_{<#2}\Par[#1]{#3}}
}}

\newcommand{\TripleRes}[4]{\MathX{%
	\Res{#1}{#4},\Res{#2}{#4},\Res{#3}{#4}%
}}

\newcommand{\Her}[1]{\MathX{%
	H(#1)%
}}

\newcommand{\Vh}[1]{\MathX{%
	V_{#1}%
}}

\newcommand{\VInter}[2]{%
	\IfEqual{#2}{\kappa}%
		{#1}
		{#1 \Inter \Vh{#2}}%
}

\newcommand{\VStr}[3][]{\MathX{%
	\Structure[#1]{\Vh{#2}, \in, #3}%
}}

\providecommand{\M}{\MathX{\Str{M}}}
\providecommand{\N}{\MathX{\Str{N}}}

\newcommand{\Model}[1]{\MathX{%
	{\M}_{#1}%
}}
\newcommand{\Voc}[1]{\MathX{%
	{\rho}_{#1}%
}}
\newcommand{\VocRes}[2]{\MathX{%
	{#1}^{{\rho}_{#2}}%
}}
\newcommand{\ModelRes}[3]{\MathX{%
	\Res{{\M}_{#1}^{\Voc{#3}}}{#2}%
}}
\newcommand{\Int}[2]{\MathX{%
	{#1}^{#2}%
}}
\newcommand{\No}[2][]{\MathX{%
	\IfEmpty{#1}%
		{\MathSymb{No}(#2)}%
		{\MathSymb{No}_{#1}(#2)}%
}}

\newcommand{\Seqs}[3][]{\MathX{%
	{\Brackets[#1]{#2}}^{#3}%
}}
\newcommand{\BSeqs}[2][]{\MathX{%
	\IfEmpty{#1}%
		{{\Brackets[]{#2}}^{<#2}_{\scriptscriptstyle B}}
		{{\Brackets[]{#2}}^{#1}_{\scriptscriptstyle B}}
}}

\newcommand{\BSeq}[2][]{\MathX{%
	\IfEmpty{#1}%
		{\boldsymbol{#2}}%
		{\boldsymbol{#2}_{#1}}%
}}

\newcommand{\Relation}[2][]{%
	\IfEmpty{#1}%
		{R_{#2}}%
		{R_{#2}^{#1}}%
}

\newcommand{\Rel}[2][]{\MathX{%
	\Relation[#1]{\BSeq{#2}}%
}}

\newcommand{\Zero}[1][]{\MathX{%
	\IfEmpty{#1}%
		{\BSeq[]{0}}%
		{\Res{\BSeq[]{0}}{#1}}%
}}

\newcommand{\T}{\MathX{T}}

\newcommand{\TLevel}[2][]{\Proj[#1]{T}{#2}}

\newcommand{\TBelow}      {\triangleleft}

\newcommand{\TBelowEq}    {\trianglelefteq}

\newcommand{\Ord}[2][]{\MathX{%
	\Symb{ord}\Par[#1]{#2}%
}}
\newcommand{\Fst}[2][]{\MathX{%
	\Symb{fun}_{1}\Par[#1]{#2}%
}}
\newcommand{\Snd}[2][]{\MathX{%
	\Symb{fun}_{2}\Par[#1]{#2}%
}}
\newcommand{\Trd}[2][]{\MathX{%
	\Symb{fun}_{3}\Par[#1]{#2}%
}}
\newcommand{\TCub}[2][]{\MathX{%
	\Symb{cst}\Par[#1]{#2}%
}}
\newcommand{\EmptyTuple}{\MathX{%
	\Tuple{\emptyset, \emptyset, \emptyset, \emptyset}%
}}

\newcommand{\TFuncs} [1] {\Functions {({\beth}_{#1})} 2}
\newcommand{\UFuncs} [1] {\Symb{Con}\Par[\big]{\TFuncs{#1}}}

\newcommand{\TLo}[1][]{%
	\IfEmpty{#1}%
		{\mathrel{\lessdot}}%
		{\mathrel{{\lessdot}_{#1}}}%
}

\newcommand{\GInds}{\MathX{%
	\Symb{Suc}^+%
}}
\newcommand{\ReflG}{\MathX{%
	\Symb{Ref}^{2}%
}}
\newcommand{\ReflH}{\MathX{%
	\Symb{Ref}^{1}%
}}
\newcommand{\Ref}[2][]{\MathX{%
	\Symb{ref}\Par[#1]{#2}%
}}
\newcommand{\Cub}[3]{\MathX{%
	C_{#1,#2,#3}%
}}

\newcommand{\U}{U}

\renewcommand{\TH}{\MathX{%
	{\T}^1%
}}
\newcommand{\UH} {\MathX{%
	{\U}^1%
}}
\newcommand{\ULevel}[2][]{\Proj[#1]{\U}{#2}}

\newcommand{\TG}[1][\lambda]{\MathX{%
	{\T}^{2}_{#1}%
}}
\newcommand{\UG}[1][\lambda]{\MathX{%
	{\U}^{2}_{#1}%
}}

\newcommand{\UGLevel} [2][\lambda] {\Proj[]{{\U}^{2}_{#1}}{#2}}

\newcommand {\UWo} {\ll}

\newcommand{\UPred}[2][]{\MathX{%
	\Proj {\U} {\UWo #2}%
}}
\newcommand{\Id}[2][]{\MathX{%
	\Symb{id}\Par[#1]{#2}%
}}
\newcommand{\EFunc}[2][]{\MathX{%
	\IfEmpty{#1}%
		{{c}^{#2}}%
		{{c}^{#2}_{#1}}%
}}

\newcommand{\E}{\MathX{\mathcal{E}}}
\newcommand{\D}{\MathX{D}}

\newcommand{\Ind}[2][]{\MathX{%
	\Symb{ind}\Par[#1]{#2}%
}}
\newcommand{\FSeq}[1]{\MathX{%
	\bar{#1}%
}}

\newcommand{\InvFSequence}[3]{\MathX{%
	\SimpleSeq{\Inv{{#1}_{#3}},\dots,\Inv{{#1}_{#2}}}%
}}

\newcommand{\Composition}[3]{\MathX{%
	{#1}_{#3} \Comp \dots \Comp {#1}_{#2}%
}}
\newcommand{\InvComposition}[3]{\MathX{%
	\Inv{{#1}_{#2}} \Comp \dots \Comp \Inv{{#1}_{#3}}%
}}
\newcommand{\ResComp}[2]{\MathX{%
	{#1}_{#2}%
}}

\newcommand{\FComps}[1][]{\MathX{%
	\IfEmpty{#1}%
		{\bar{\mathcal{F}}}
		{\mathcal{F}_{#1}}
}}

\newcommand{\F}{\FComps[1]}

\newcommand{\FuncSeq}[3][]{\MathX{%
	\IfEmpty{#1}%
		{\FSeq{g}^{\bar{#2},#3}}%
		{{g}^{\bar{#2},#3}_{#1}}%
}}
\newcommand{\Func}[3][]{\MathX{%
	\IfEmpty{#1}%
		{{g}^{\bar{#2},#3}}%
		{{g}^{\bar{#2},#3}_{#1}}%
}}

\newcommand{\FCompSeqs}{\MathSymbX{Seq}}

\newcommand{\SLh}[2][]{\MathX{%
	\Symb{lh}(\bar{#2}#1)%
}}
\newcommand{\BLh}[2][]{\MathX{%
	\Symb{lh}(\BSeq{#2}#1)%
}}

\newcommand{\Beg}[2][]{\MathX{%
	\Symb{beg}\Par[#1]{#2}%
}}
\newcommand{\End}[2][]{\MathX{%
	\Symb{end}\Par[#1]{#2}%
}}
\begin{document}

%

\title{\protect
The number of $L_{\infty\kappa}$-equivalent %
non-isomorphic models for $\kappa$ weakly compact%
}

\author{%
	Saharon Shelah%
	\thanks{%
	Research supported by the United States-Israel Binational
	Science Foundation. Publication 718.
	}\\%
	\and%
	Pauli V\"{a}is\"{a}nen%
	\thanks{%
		The second author wishes to thank Tapani Hyttinen
		under whose supervision he did his share of the paper.
	}%
}

\date{November 12, 1999}

\maketitle

\begin{abstract}%
For a cardinal $\kappa$ and a model \M of cardinality $\kappa$ let \No
\M denote the number of non-isomorphic models of cardinality $\kappa$
which are \Lan \kappa-equivalent to \M.  In \cite {Sh133} Shelah
established that when $\kappa$ is a weakly compact cardinal and $\mu
\leq \kappa$ is a nonzero cardinal, there exists a model \M of
cardinality $\kappa$ with $\No \M = \mu$.  We prove here that if
$\kappa$ is a weakly compact cardinal, the question of the possible
values of \No \M for models \M of cardinality $\kappa$ is equivalent
to the question of the possible numbers of equivalence classes of
equivalence relations which are $\Sigma^1_1$-definable over \Vh
\kappa. In \cite {ShVaSigma} we proved that, consistent wise, the
possible numbers of equivalence classes of $\Sigma^1_1$-equivalence
relations can be completely controlled under the singular cardinal
hypothesis. These results settle the problem of the possible values of
\No \M for models of weakly compact cardinality, provided that the
singular cardinal hypothesis holds.
 \footnote{%
 1991 Mathematics Subject Classification: %
	primary 03C55; secondary 03C75. %
Key words: %
	number of models, infinitary logic.
}%
 \end{abstract}
%

\begin{SECTION} {-} {Introduction} {Introduction}
%
Suppose $\kappa$ is a cardinal and \M is a model of cardinality
$\kappa$. Let \No \M denote the number of non-isomorphic models of
cardinality $\kappa$ which are elementary equivalent to \M over the
infinitary language \Lan \kappa. We study the possible values of \No
\M for different models \M.

%
When \M is countable, $\No \M = 1$ by \cite {Scott}. This result
extends to all structures of singular cardinality $\lambda$ provided
that $\lambda$ is of countable cofinality \cite {Chang}.  The case \M
is of singular cardinality $\lambda$ with uncountable cofinality
$\kappa$ was first treated in \cite {Sh189} and later on in \cite
{Sh228}. In these papers Shelah showed that if $\kappa > \aleph_0$,
$\theta^\kappa < \lambda$ for every $\theta < \lambda$, and $0 < \mu <
\lambda$ or $\mu = \lambda^\kappa$, then $\No \M = \mu$ for some model
\M of cardinality $\lambda$. In the paper \cite {ShVa644} of the
authors the singular case is revisited, and particularly, it is
established, under the same assumptions as above, that the values
$\mu$ with $\lambda \leq \mu < \lambda^\kappa$ are possible for \No \M
with \M of cardinality $\lambda$.

If $V = L$, $\kappa$ is an uncountable regular cardinal which is not
weakly compact, and \M is a model of cardinality $\kappa$, then $\No
\M \in \Braces {1, 2^\kappa}$ \cite {Sh129}. For $\kappa = \aleph_1$
this result was first proved in \cite {Palyutin}. The values $\No \M
\in \Braces {\aleph_0, \aleph_1}$ for a model of cardinality $\aleph
_{1}$ are consistent with $\ZFC + \GCH$ as noted in \cite {Sh129}. All
the nonzero finite values of \No \M for models of cardinality $\aleph
_{1}$ are proved to be consistent with $\ZFC + \GCH$ in \cite
{ShVa646}.


When $\kappa$ is a weakly compact cardinal and $\mu$ is a nonzero
cardinal $\leq \kappa$ there is a model \M of cardinality $\kappa$
with $\No \M = \mu$ \cite {Sh133}. In the present paper we show that
when $\kappa$ is a weakly compact cardinal, the possible values of \No
\M for models of cardinality $\kappa$ depends only on the possible
numbers of equivalence classes of equivalence relations which are
$\Sigma^1_1$-definable over \Vh \kappa as follows: for some first
order sentence $\phi$ in the vocabulary \Braces {\in, R_0, R_1, R_2,
R_3} and a subset \Param of \Vh \kappa, it is the case that for all
$s, t \in \Functions \kappa 2$,
 \[
	s \Equiv[] t
	\Text{iff for some} r \in \Functions \kappa 2 \ 
	\VStr \kappa {\Param, s, t, r} \models \phi
 \]
where \Param, $r$, $s$, and $t$ are the interpretations of the symbols
$R_0$, $R_1$, $R_2$, and $R_3$ respectively. More precisely, we prove
the following theorem.

\begin{THEOREM}{NoMiffSigma^1_1}%
When $\kappa$ is a weakly compact cardinal, the following two
conditions are equivalent for every nonzero cardinal $\mu$:%
\begin{ITEMS}%

\ITEM{NoEquiv}%
there is an equivalence relation on \Functions \kappa 2 which is
$\Sigma_1^1$-definable over \Vh \kappa and has exactly $\mu$ different
equivalence classes;

\ITEM{NoM}%
$\No \M = \mu$ for some model \M of cardinality $\kappa$.

\end{ITEMS}%
\end{THEOREM}

In the paper \cite {ShVaSigma} we proved that for every nonzero
cardinal $\mu \in \kappa \Union \Braces {\kappa, \SuccCard \kappa,
2^\kappa}$ there is always a $\Sigma^1_1$-equivalence relation \Note
{as defined above} with exactly $\mu$ different equivalence
classes. Moreover, consistent wise, one can completely control the
possible numbers of equivalence classes of $\Sigma^1_1$-equivalence
relations provided that the singular cardinal hypothesis holds \cite
[Theorem 1] {ShVaSigma}. It follows that, the question of possible
value of \No \M is completely solved, when \M is of weakly compact
cardinality and the singular cardinal hypothesis holds. Again more
formally, the conclusion will be the following.

\begin{CONCLUSION}{NoM}%
Suppose that the following conditions are satisfied:
\begin{itemize}%

\item $\kappa$ is a weakly compact cardinal;

\item $\kappa$ remains a weakly compact cardinal in the standard Cohen
forcing adding a new subset of $\kappa$;

\item the singular cardinal hypothesis holds;

\item $\lambda > \SuccCard \kappa$ is a cardinal with $\lambda^\kappa =
\lambda$;

\item $\Omega$ is a set of cardinals between \SuccCard \kappa and
$\lambda$ \Note {possibly empty}, which is closed under unions of
$\leq \kappa$-many cardinals and products of $< \kappa$-many
cardinals.

\end{itemize}%
Then there is a forcing extension where there are no new sets of
cardinality $< \kappa$, all cardinals and cofinalities are preserved,
$\kappa$ remains a weakly compact cardinal, $2^\kappa = \lambda$, and
for all cardinals $\mu$, there exists a model \M of cardinality
$\kappa$ with $\No \M = \mu$ if and only if $\mu$ is a nonzero
cardinal $\leq \SuccCard \kappa$ or $\mu$ is in $\Omega \Union
{\Braces {2^\kappa}}$.
\end{CONCLUSION}

\Remark When $\kappa$ is a weakly compact cardinal, it is possible to
have, using the upward Easton forcing, a generic extension where
$\kappa$ is still a weakly compact cardinal and $\kappa$ remains
weakly compact in the Cohen forcing adding a subset of $\kappa$
(Silver). The forcing needed in the conclusion is the ordinary way to
add Kurepa trees of height $\kappa$ with $\mu$-many $\kappa$-branches
through them, for all $\mu \in \Omega$. As noted in \cite [Fact 5.1]
{ShVaSigma}, this forcing is locally $\kappa$ Cohen, and therefore,
$\kappa$ remains a weakly compact cardinal in the composite forcing of
the upward Easton forcing and the addition of new Kurepa trees.

Note also, that the closure properties mentioned in the conclusion are
necessary by the fact that the possible numbers of equivalence classes
of $\Sigma_1^1$-equivalence relations are always closed under unions
of length $\leq \kappa$ and products of length $< \kappa$, see \cite
[Lemma 3.4] {ShVaSigma}.

There are three parts in the paper. First in \Section {Preliminaries}
we recall a definition of a \EhrFra-game \EF \kappa \lambda \M \N
generalizing the elementary equivalence between two models over an
infinitary language \Lan \kappa.  In addition to that, we show that if
\M is a model of cardinality $\kappa$, then there is a
$\Sigma^1_1$-equivalence relation with exactly \No \M different
equivalence classes. We also note how this connection extends to even
more strongly notions of equivalence between models.

The last two sections are dedicated to the other half of the proof of
the theorem, namely to the proof that the existence of a
$\Sigma^1_1$-equivalence relation with $\mu$-many equivalence classes
implies the existence of a model \M with $\Card \M = \kappa$ and $\No
\M = \mu$ \Note {\Lemma {No} at the end of \Section {Models}}. First
in \Section {Functions} we define a special family of functions which
is used to build models in the last section. This part might feel
quite technical. However, it is elementary and the fundamental idea of
the construction is from \cite {Sh133}. The reader may even skip all
the lemmas of this section in the first reading, and return to those
when they are referred from the last section.

The content of \Section {Models} is as follows.  Assuming that
$\kappa$ is strongly inaccessible and $\phi$ defines a
$\Sigma^1_1$-equivalence relation $\Equiv$ on \Functions \kappa 2 with
a parameter $\Param \Subset \Vh \kappa$ we construct models \Model t
for $t \in \Functions \kappa 2$ satisfying that
\begin{itemize}

\item%
	the models are of cardinality $\kappa$ and they have a common
	vocabulary $\rho$ consisting of $\kappa$ many relation symbols
	each of having arity $< \kappa$;

\item%
	all the models are pairwise \Lan \kappa-equivalent, and even
	more, they are pairwise \MLan \kappa \lambda-equivalent for
	any previously fixed regular cardinal $\lambda < \kappa$ \Note
	{see \Definition {EF}};

\item%
	for all $s, t \in \Functions \kappa 2$, the models \Model s
	and \Model t are isomorphic if, and only if, $s$
	and $t$ are equivalent with respect to $\Equiv$.

\end{itemize}
Furthermore, when $\kappa$ is a weakly compact cardinal the models
satisfy the additional property that
\begin{quote}

	if a model \N has vocabulary $\rho$, \N is of cardinality
	$\kappa$, and \N is \Lan \kappa-equivalent to one (all) of the
	models \Model t, $t \in \Functions \kappa 2$, then \N is
	isomorphic to \Model s for some $s \in \Functions \kappa 2$.

\end{quote}
This is the main difference between the strongly inaccessible
non-weakly compact case and the weakly compact case: the
$\Pi^1_1$-indescribability property of a weakly compact cardinal
$\kappa$ ensures that the ``isomorphism type'' of any model \N with
domain $\kappa$ is already determined by the isomorphism types of the
bounded parts \Res \N \alpha, $\alpha < \kappa$, alone \Note {see the
proof of \Lemma {No(M)<2^k}}.
%
\end{SECTION}

\begin{SECTION} {-} {Preliminaries} {Preliminaries}
%


\begin{DEFINITION}{EF}%
Suppose $\mu$ is a cardinal and $\lambda$ is a infinite regular
cardinal. Let \M and \N be models of a common relational
vocabulary. The \EhrFra-game \EF \mu \lambda \M \N is defined as
follows.  The game has two players, \PlayerOne and \PlayerTwo, and a
play of the game continues for at most $\lambda$ rounds. On round $i <
\lambda$ player \PlayerOne first chooses $X_i \in \Braces {\M, \N}$
and $A_i \Subset X_i$ of cardinality $< \mu$. Then \PlayerTwo replies
with a partial isomorphism $p_i$ such that
\begin{itemize}
\item $\Dom {p_i} \Subset \M$, $\Ran {p_i} \Subset \N$, $\BigUnion [j < i]
{p_j} \Subset p_i$, and

\item $A_i \Subset \Dom {p_i}$ if $X_i = \M$, and $A_i \Subset \Ran
{p_i}$ otherwise.

\end{itemize}
Player \PlayerTwo wins a play if the play lasts $\lambda$ many rounds,
and otherwise, \PlayerOne wins the play. We write $\M \MEquiv \kappa
\lambda \N$ when \PlayerTwo has a winning strategy in \EF \kappa
\lambda \M \N.
\end{DEFINITION}

Let \M and \N be models of a common relational vocabulary and $\kappa$
be a cardinal. The game \EF \kappa \omega \M \N is the usual
\EhrFra-game of length $\omega$ which characterizes the existence of a
nonempty family of partial isomorphism with the ``fewer than $\kappa$
at the time back-and-forth property''. If \M and \N satisfy the same
sentences of the infinitary language \Lan \kappa, we write $\M \LEquiv
\kappa \N$. By the Karp's theorem \cite {Karp65} player \PlayerTwo has
a winning strategy in \EF \kappa \omega \M \N if, and only if, $\M
\LEquiv \kappa \N$. So the game \EF \kappa \lambda \M \N, for an
infinite regular cardinal $\lambda < \kappa$, is a generalized version
of the ``fewer than $\kappa$ at the time back-and-forth property''.
There are so-called infinitely deep languages \MLan \kappa \lambda
with the property that $\M \MEquiv \kappa \lambda \N$ if, and only if,
\M and \N satisfy the same sentences of \MLan \kappa \lambda \cite
{Hytt90,Kartt,Oikk97}.

For a model \N we let \Card \N denote the cardinality of the universe
of \N.  For any model \M of cardinality $\kappa$ and a regular
cardinal $\aleph_0 \leq \lambda < \kappa$, we define $\No [\lambda]
\M$ to be the cardinality of the set
 \[
	\Set [\big] {\Quotient \N \Isomorphic}
	{ \Card \N = \kappa \And \N \MEquiv \kappa \lambda \M },
 \]
where \Quotient \N \cong is the equivalence class of \N under the
isomorphism relation.

For all sets $X$ of ordinals, the ordinal $\sup \Set {\alpha+1}
{\alpha \in X}$ is abbreviated by \PSup X. For all sequences $\bar
\alpha = \Seq {\alpha_i} {i < \theta}$ of ordinals, we denote $\PSup
[\big] {\Set {\alpha_i} {i < \theta}}$ by $\PSup {\bar \alpha}$, and
we abbreviate the sequence \Seq {f(\alpha_i)} {i < \theta} by $f(\bar
\alpha)$. For a regular cardinal $\kappa$ and a subset $S$ of
$\kappa$, $S$ is called stationary if for every closed and unbounded
subset $C$ of $\kappa$, $S \Inter C$ is nonempty.


Next we recall the definition of an equivalence relation which is
$\Sigma^1_1$-definable over the set \Her \kappa of all sets
hereditarily of cardinality $< \kappa$ from \cite [Definition 3.1]
{ShVaSigma}. In this paper $\kappa$ will be a strongly inaccessible
cardinal and so \Her \kappa equals to the set \Vh \kappa of all sets
having rank $< \kappa$. It will be more convenient to use \Vh \kappa
instead of \Her \kappa here, especially, when we consider elementary
submodels of \Vh \kappa, and also, when we apply
$\Pi^1_1$-indescribability property for $\kappa$ weakly compact.

\begin{DEFINITION}{Sigma}
Suppose $\kappa$ is a strongly inaccessible cardinal. We say that
$\phi$ defines a $\Sigma^1_1$-equivalence relation $\Equiv$ on
\Functions \kappa 2 with a parameter $P \Subset \Vh \kappa$ when
\begin{ITEMS}

\ITEM{voc}%
$\phi$ is a first order sentence in a vocabulary consisting of $\in$,
one unary relation symbol $R_0$, and binary relation symbols $R_1$,
$R_2$, and $R_3$;

\ITEM{eq}%
the following definition gives an equivalence relation on \Functions
\kappa 2: for all $s,t \in \Functions \kappa 2$ %
 \[
	s \Equiv t \Iff \ 
	\VStr \kappa {P, s, t, r} \models \phi
	\ForSome r \in \Functions \kappa 2,
 \]
where $P$, $s$, $t$, and $r$ are the interpretations of the symbols
$R_0$, $R_1$, $R_2$, and $R_3$ respectively.

\end{ITEMS}

We say that there exists a $\Sigma^1_1$-equivalence relation on
\Functions \kappa 2 having $\mu$ many different equivalent classes
when there is some sentence $\phi$ and a parameter $P \Subset \Vh
\kappa$ such that $\phi$ defines a $\Sigma^1_1$-equivalence relation
$\Equiv$ on \Functions \kappa 2 with the parameter $P$ and $\Card
[\big] {\Set {\Quotient f {\Equiv}} {f \in \Functions \kappa 2}} =
\mu$.
\end{DEFINITION}

\begin{LEMMA}{Sigma}%
If $\kappa$ is a strongly inaccessible cardinal and \M is a model of
cardinality $\kappa$, then there is a $\Sigma_1^1$-equivalence
relation on \Functions \kappa 2 such that the number of different
equivalence classes of it is $\No \M$.

\begin{Proof}%
For a function $\pi$ having domain $\kappa$ and a binary relation $R$,
we let $\pi(R)$ denote the set \Set {\pi(\xi)} {\ForSome \xi < \kappa,
\Pair \xi 1 \in R}. For every $n < \omega$, we write $\pi_n$ for a
fixed definable bijection from $\kappa$ onto \Set {\Tuple{\alpha_1,
\dots, \alpha_n}} {\alpha_1, \dots, \alpha_n \in \kappa}. We may
assume that the domain of \M is $\kappa$ and its vocabulary consists
of one relation symbol $Q$ of finite arity $n$. For a binary relation
$R$ let $\N(R)$ be the model having domain $\kappa$ and interpretation
$\pi_n(R)$ for the relation symbol $Q$. By the inaccessibility of
$\kappa$, let $\rho$ be a bijection from $\kappa$ onto \Vh \kappa. For
a binary relation $R$ let $\tau_1(R)$ be the set \Set {\rho(\xi)} {\xi
< \kappa \Text {is a successor ordinal and} \Pair \xi 1 \in R} and
$\tau_2(R)$ be the set \Set {\rho(\xi)} {\xi < \kappa \Text {is a
limit ordinal and} \Pair \xi 1 \in R}.

Since for all models \N, the game \EF \kappa \omega \M \N is
determined, the condition $\M \not\LEquiv \kappa \N$ is equivalent to
that player \PlayerOne has a winning strategy in \EF \kappa \omega \M
\N. Therefore, using the interpretation $Q^\M$ and the bijection
$\rho$ as a parameter, the sentence $\phi(R_0,R_1, R_2, R_3)$ saying
\begin{quote}%
``\big (
$\tau_1(R_3)$ is a winning strategy for player \PlayerOne in \EF
\kappa \omega \M {\N(R_1)} and
$\tau_2(R_3)$ is a winning strategy for player \PlayerOne in \EF
\kappa \omega \M {\N(R_2)}
\big )
or
\big (
$\pi_2(R_3)$ is an isomorphism between $\N(R_1)$ and $\N(R_2)$
\big )''
\end{quote}%
defines a $\Sigma^1_1$-equivalence relation on \Functions \kappa 2.
This definition is as wanted in the claim, except that when \No \M is
finite the definition gives one extra class. That can be avoided by
obvious changes to the definition.
\end{Proof}
\end{LEMMA}

\Remark As noted in \cite [Section 5] {ShVaSigma}, the theorem on the
possible numbers of equivalence classes of $\Sigma^1_1$-equivalence
relations directly extends to equivalence relations which are
definable over \Vh \kappa using a subset of \Vh \kappa as a parameter
and a sentence which is a Boolean combination of a sentence containing
one second order existential quantifier \Note {$\Sigma^1_1$-sentence}
and a sentence containing one second order universal quantifier \Note
{$\Pi^1_1$-sentence}.
In the proof of \Lemma {Sigma} above we needed the fact that the game
\EF \kappa \omega \M \N is determined to find a $\Sigma^1_1$-sentence
which says ``if the models are equivalent then \dots'', i.e., ``either
the models are nonequivalent or \dots''.
But there is a $\Pi^1_1$-sentence saying ``there is no winning
strategy for player \PlayerTwo in the game \dots''. Hence, after we
have proved that existence of a $\Sigma^1_1$-equivalence relation with
$\mu$-classes implies existence of a model \M with $\No \M = \mu$
\Note {after the next two sections}, we can conclude: consistent wise,
the possible values of \No \M for \M of weakly compact cardinality
$\kappa$ might be exactly as wanted, and moreover, the possible values
of \No [\lambda] \M, for \M of cardinality $\kappa$ and $\lambda <
\kappa$ any regular cardinal, coincide with the possible values of \No
[\omega] \M.
%
\end{SECTION}

\begin{SECTION} {-} {Functions} {The family of functions}
%

Throughout the next two sections $\kappa$ is a strongly inaccessible
cardinal, i.e., a regular limit cardinal satisfying $2^\mu < \kappa$
for all $\mu < \kappa$, and $\lambda$ is a fixed regular cardinal
below $\kappa$. %
\PvComment {$\lambda$ is first time used in \Definition {G}.}
\PvComment {$\kappa$ is weakly compact in \Lemma {No(M)<2^k}}

In this section we define a family of functions which will be used to
build the models \Model t, $t \in \Functions \kappa 2$ \Note
{\Definition {Models}}. There is a similar idea in \cite {Sh133},
however, this time we want the functions to satisfy some additional
properties. Hence the definition of the family will be a little bit
more complicated. Particularly, we shall first define a special tree
\Note {\Definition {U}} which will be only a steering apparatus in the
construction of the family of functions itself \Note {\Definition
{p}}.

To make our models strongly equivalent we shall guarantee that for
every pair \Model s and \Model t, $s,t \in \Functions \kappa 2$, a
certain subfamily of all the functions will form a winning strategy
for player \PlayerTwo in the game \EF \kappa \lambda {\Model s}
{\Model t} \Note {\Definition {G}}. So we shall need a ``stronger
extension property'' for the functions than was needed in \cite
{Sh133}.

Most importantly, we want that a pair \Model s and \Model t, $s,t \in
\Functions \kappa 2$, of models are isomorphic if, and only if, the
corresponding indices $s$ and $t$ are equivalent with respect to some
previously fixed $\Sigma^1_1$-equivalence relation \Note {\Lemma
{Characterization}}. Hence we have to code information about the
equivalence relation into the family of functions \Note {\Definition
{H} and \Definition {p}}.

Henceforth $\phi$ denotes a sentence which defines a
$\Sigma^1_1$-equivalence relation $\Equiv$ on \Functions \kappa 2 with
a parameter $P \Subset \Vh \kappa$ \Note {see \Definition
{Sigma}}. Without loss of generality we may assume that for all $s,t,r
\in \Functions \kappa 2$,
 \EQUATION{symmetry}{
	\VStr \kappa {P, s, t, r} \models \phi \Iff
	\VStr \kappa {P, t, s, r} \models \phi.
 }%
\SimpleComment{%
Namely, if $\psi$ defines a $\Sigma^1_1$-equivalence relation on
\Functions \kappa 2 with a parameter $P$ but does not have the
property wanted, then define $\phi$ to be the formula $(\psi \Or
\psi')$ where $\psi'$ is the same sentence as $\psi$ except that $R_1$
is replaced with $R_2$ and vice versa. Then $\phi$ defines the same
$\Sigma^1_1$-equivalence relation on \Functions \kappa 2 as $\psi$
defines, and additionally, $\phi$ satisfies \Equation {symmetry}.
}
Furthermore, we may assume that for all $s,t \in \Functions
\kappa 2$,
 \EQUATION{constant}{
	\If s(\beth_\alpha) = t(\beth_\alpha) \ForEvery \alpha < \kappa
	\Then s \Equiv t.
 }%
\SimpleComment{%
Again, if $\psi$ defines a $\Sigma^1_1$-equivalence relation $\Equiv
[\psi, P]$ on \Functions \kappa 2 with a parameter $P$ but does not
have the property wanted, then define $\phi$ to be the sentence where
each occurrence of ``$\Pair \alpha i \in R_l$'' in $\psi$, for $\alpha
< \kappa$, $i \in \Braces {0, 1}$, and $l \in \Braces {1, 2, 3}$, is
replaced with ``$\Pair {\beth_\alpha} i \in R_l$''. Then $\Equiv
[\phi, P]$ and $\Equiv [\psi, P]$ has the same number of equivalence
classes and $\phi$ satisfies \Equation {constant}.
}


%
\begin{DEFINITION}{T}%
Let \TLevel 0 be \Braces \EmptyTuple and for every nonzero $\alpha <
\kappa$ define
 \[
	\TLevel \alpha = \Set {\Tuple{\eta,\nu,\tau,C}}
	{\eta, \nu, \tau \in \TFuncs \alpha, \eta \not= \nu, 
	\And C \Text{is a closed subset of} \alpha},
 \]
$\TLevel {<\alpha} = \BigUnion [\beta < \alpha] {\TLevel \beta}$, and
$\T = \BigUnion [\alpha < \kappa] {\TLevel \alpha}$.
\PvComment {$C$ ensures that all chains of length $< \lambda$ will be
extended, see \Definition {G}.}

For every $u \in \T$, we let \Ord u, \Fst u, \Snd u, \Trd u, and \TCub
u be elements such that $u \in \TLevel {\Ord u}$ and $u = \Tuple {\Fst
u, \Snd u, \Trd u, \TCub u}$. Furthermore, for every $u \in \T$,
we let \Res u \beta denote \EmptyTuple when $\beta = 0$, and when
$\beta > 0$,
\[
	\Res u \beta = \Tuple {\Res {\Fst u} {\beth_\beta},
	\Res {\Snd u} {\beth_\beta},
	\Res {\Trd u} {\beth_\beta},
	\TCub u \Inter \beta}.
\]
The elements $u,v \in \T$ form a tree when they are ordered by
 \[
	u \TBelow v \Iff
	u = \Res v {\Ord u} \And \Ord u \in \TCub v.
 \]
The notation $u \TBelowEq v$ stands for $u \TBelow v$ or $u = v$. For
a $\TBelow$-increasing chain \Seq {u_i} {i < \theta} of elements in
\T, $\theta < \kappa$, we write \BigUnion [i < \theta] {u_i} for the
following element in \T:
 \ARRAY{
	\Tuple [\Big] {
		\BigUnion [i < \theta] {\Fst {u_i}},
		\BigUnion [i < \theta] {\Snd {u_i}},
		\BigUnion [i < \theta] {\Trd {u_i}},
		C
	},
 }
where $C$ is the closure of \BigUnion [i < \theta] {\TCub {u_i}}.
\PvComment{union $+$ sup if sub C $<$ sup sequence}%
\end{DEFINITION}

\begin{DEFINITION}{H}
For every $s, t, r \in \Functions \kappa 2$ define
 \[
	\Cub s t r = \Braces 0 \Union
	\Set {\delta < \kappa}
	{\VStr \delta {\VInter P \delta, \TripleRes str \delta}
	\ElemSubstr
	\VStr \kappa {P,s,t,r} \models \phi }.
 \]
\Note {Then for all nonzero $\delta \in \Cub s t r$, $\beth_\delta =
\delta$.}  We let \TH be the set of all $u \in \T$ such that for some
$s, t, r \in \Functions \kappa 2$, the following conditions are
satisfied:
\begin{itemize}

\item	$\Fst u \Subset s$,
	$\Snd u \Subset t$, and
	$\Trd u \Subset r$;

\item	$\Ord u \in \Cub s t r$ and
	$\TCub u = \Cub s t r \Inter \Ord u$.

\end{itemize}
\PvComment{We need $\EmptyTuple \TBelow u$}%
\end{DEFINITION}

\begin{DEFINITION}{TLo}
For each $\alpha < \kappa$ define a lexicographic order $\TLo
[\alpha]$ as follows: for all elements $\eta, \nu \in \TFuncs \alpha$,
 \[
	\eta \TLo [\alpha] \nu \Iff
	\eta \not= \nu \And
	\eta(\xi) < \nu(\xi), \For
	\xi = \min
	\Set {\zeta < \beth_\alpha} {\eta(\zeta) \not= \nu(\zeta)}.
 \]
Define
 \ARRAY{
	\Proj \T \TLo = \Set {u \in \T} {\Fst u \TLo [\Ord u] \Snd u}.
 }

For every $u \in \T$, denote the tuple \Tuple {\Snd u, \Fst u, \Trd u,
\TCub u} by \Inv u \Note {the order of the first and the second
elements are exchanged}. %
\PvComment{These definitions are used only for notational
convenience. In particular, because of \Definition {FComps}.}
\end{DEFINITION}

\Remark For every $u \in \TH$, $\Inv u \in \TH$ by the assumption
\Equation {symmetry}.

\begin{DEFINITION}{G}
Define \GInds to be the set \Set {\beta+1} {\beta \Text {is a
successor ordinal\,}}. %
\PvComment{%
Note that \GInds is disjoint from the set \Set {\Ord u} {u \in
\TH}.
}
Assume $\pi$ is a surjective function from $\Braces 0 \Union \GInds$
onto $\Braces \EmptyTuple \Union \Proj \T \TLo$ such that
\begin{itemize}

\item if $\pi(\alpha) = u$ then either $\alpha=0$ and $u =
\EmptyTuple$, or else $\Ord u < \alpha$;

\item for every $u \in \Braces \EmptyTuple \Union \Proj \T \TLo$, the
set \Set {\alpha \in \GInds} {\pi(\alpha) = u} is unbounded in
$\kappa$.

\end{itemize}
We define \TG, for fixed regular $\lambda < \kappa$, to be the
smallest subset of \T satisfying the following conditions.

\begin{enumerate}

\ITEM{empty} \EmptyTuple is in \TG.

\ITEM{inv} If $u \in \TG$ then $\Inv u \in \TG$.

\ITEM{other}%
\TG contains every $u \in \Proj \T \TLo$ having the properties:
\begin{enumerate}

\ITEM{successor}%
If $\sup \TCub u < \Ord u$ then $\Ord u \in \GInds$ and for the
maximal element $\gamma \in \TCub u$, $\pi \Par[\big] {\Ord
u} = \Res u \gamma \in \TH \Union \TG$;
\PvComment {$u \in \Proj \T \TLo$ implies $\Res u \beta \in \Proj \T
\TLo$}%
\PvComment {The elements of \TH also need to be extended.}

\ITEM{limit}%
if $\sup \TCub u = \Ord u$, then $\TCub u \Inter \GInds$ is nonempty,
$\Cf {\Ord u} < \lambda$, and for every $\beta \in \TCub u$, $\Res u
\beta \in \TH \Union \TG$.

\end{enumerate}
\end{enumerate}
\end{DEFINITION}

We shall need only a restricted part of \T, so we change the notation
a little bit.

\begin{DEFINITION}{U}%
For each $\alpha < \kappa$, let \UFuncs \alpha be the family of
functions $\eta$ satisfying that $\eta \in \TFuncs \alpha$ for some
$\alpha < \kappa$, $\eta(\xi) = 0$ if $\xi < \beth_0$, and for all
$\beta < \alpha$ and $\xi \in \beth_{\beta+1} \Minus \beth_{\beta}$,
$\eta(\xi) = \eta(\beth_\beta)$. Define the restricted part of \T to
be
\ARRAY[ll]{
	\U =&
	\Set [\Big] {u \in \TH \Union \TG}
	{u = \EmptyTuple \Or \\
	& \Fst u, \Snd u, \Trd u \in \UFuncs {\Ord u}}.
}%
Denote $\TH \Inter \U$ by \UH and $\TG \Inter \U$ by \UG.  We write
that $u \in \U$ is a successor of $v$ when $v \in \U$ and there is no
$w \in \U$ with $v \TBelow w \TBelow u$. An element $u \in \U$ is
called a successor if $u = \EmptyTuple$ or there is $v \in U$ such
that $u$ is a successor of $v$. If $u$ is not a successor, $u$ is
called a limit. For all $\alpha < \kappa$, \ULevel \alpha denote
$\TLevel \alpha \Inter \U$ and \ULevel {< \alpha} stand for $\TLevel
{< \alpha} \Inter \U$. When $u \in \U$, we write that ``\,for all $v
\TBelow u$'' when we mean that ``\,for all $v \in \U$ with $v \TBelow
u$''.
\end{DEFINITION}

\Remark By the assumption \Equation {constant} at the beginning of
this section and the use of elementary submodels in \Definition {H},
our restriction of \TFuncs \alpha to \UFuncs \alpha is
harmless. However this restriction turned out to be useful in the
forthcoming \Definition {p} and the property \ItemOfLemma {Ran} {seq}
\Note {we want that for successor $u \in \U$ with \EFunc u increasing,
the information \End {p_u} is determined by a single point $\zeta \in
\Ran {\EFunc u}$ alone, see the definitions \DefinitionNo {p} and
\DefinitionNo {FComps} below for unexplained notation}.

Note also that an element $u \in \UH$ is a successor of $v \in \U$
only if $v \in \UH$. However, for every $v \in \UH$ there is $u \in
\UG$ which is a successor of $v$. When $u \in \UH$ is a limit point
then it is a limit of elements in \UH. Besides, $u \in \UG$ is a limit
point only if it is a limit of elements in \UG %
\Note {see the proof of \ItemOfFact {U} {H-limit} below}.
\PvComment {Note, that if $u \in \UH$ and \Ord u is a regular
cardinal, then $u$ is a limit point.}

\begin{FACT}{U}
\begin{ITEMS}

\ITEM{G-closed}
For all $\TBelow$-increasing chain \Seq {u_i} {i < \theta} of elements
in \UG, $\Cf \theta < \lambda$ and the tuple \BigUnion [i < \theta]
{u_i} is in \UG.

\ITEM{H-limit}
For all $\TBelow$-increasing chain $\bar u = \Seq {u_i} {i < \theta}$
of elements in \UH, if $\bar u$ has some upper bound in \U \Note {with
respect to the order $\TBelow$} then the tuple \BigUnion [i < \theta]
{u_i} is in \UH.

\end{ITEMS}

\begin{Proof}%
\ProofOfItem{G-closed}%
This property is an obvious consequence of \ItemOfDefinition {G}
{limit}.

\ProofOfItem{H-limit}%
Suppose that $w \in \U$ is an upper bound for $\bar u$.  Let $v$ be
the smallest element in \U with $v \TBelowEq w$ and $u_i \TBelow v$
for every $i < \theta$. Then $v$ is a limit of the elements $u_i \in
\UH$, $i < \theta$.

Suppose $v$ is in \UG. By \ItemOfDefinition {G} {successor}, \TCub v
does not contain a maximal element, $\TCub v = \BigUnion [i < \theta]
{\TCub {u_i}}$ since $\bar u$ is $\TBelow$-increasing, and for every
$\beta \in \TCub u$, $\Res u \beta \in \U$.  Since \GInds contains
only successor ordinals and each $u_i$ is in \UH, $\TCub {u_i} \Inter
\GInds = \emptyset$ for every $i < \theta$.  Hence \TCub v is disjoint
from \GInds contrary to \ItemOfDefinition {G} {limit}.

It follows that $v$ must be in \UH and there are $s,t,r \in \Functions
\kappa 2$ such that
 \ARRAY[lll]{
	\Fst v \Subset s,
	\Snd v \Subset t,
	\Trd v \Subset r, \\
	\Ord v \in \Cub s t r, \And
	\TCub v = \Cub s t r \Inter \Ord v.
 }
For every $i < \theta$ and $\alpha_i = \Ord {u_i}$, $u_i \TBelow v$
implies that $\alpha_i \in \TCub v \Subset \Cub s t r$. Hence, for
each $i < \theta$,
 \[
	\VStr {\alpha_i} {\VInter P {\alpha_i}, \TripleRes str {\alpha_i}}
	\ElemSubstr \VStr \kappa {P,s,t,r} \models \phi,
 \]
and for $\delta = \BigUnion [i < \theta] {\alpha_i}$,
 \[
	\VStr \delta {\VInter P \delta, \TripleRes str \delta}
	\ElemSubstr
	\VStr \kappa {P,s,t,r} \models \phi.
 \]
Consequently, the tuple $\Tuple {\TripleRes str \delta, \Cub str
\Inter \delta} = \BigUnion [i < \theta] {u_i}$ is in \UH \Note{note
that $\beth_{\alpha_i} = \alpha_i$, $\beth_\delta = \delta$, and by
the choice of $v$, $v = \BigUnion [i < \theta] {u_i}$}.
\end{Proof}
\end{FACT}
%


%
\begin{DEFINITION}{E}%
For every $\beta < \kappa$ and $\gamma < \beth_{\beta+1}$ we define a
function $\EFunc [\gamma] \beta$ with domain $\beth_\beta$ as follows:
for all $\xi < \beth_\beta$,
 \[
	\EFunc [\gamma] \beta (\xi) =
	       \Par[\big]{\beth_\beta \Times (\gamma+1)} + \xi,
 \]
where $\Times$ and $+$ are the ordinal multiplication and addition
respectively. Write \E for the family of functions \Set {\Res f A} {f
\in E \And A \Subset \Dom f}, where $E$ is the set \BigUnion {\Set
[\big] { \Braces {\EFunc [\gamma] \beta, \Inv {(\EFunc [\gamma]
\beta)}}} {\beta < \kappa \And \gamma < \beth_{\beta +1}}}. The
reflection point of $d \in \E \Minus \Braces \emptyset$, denoted be
\Ref d, is the unique ordinal $\beta$ for which there is $\gamma <
\beth_{\beta +1}$ satisfying that either $d \Subset \EFunc [\gamma]
\beta$ or $d \Subset \Inv {\Par [\big] {\EFunc [\gamma] \beta}}$.
\end{DEFINITION}

\begin{FACT}{E}%
\begin{ITEMS}%
\ITEM {increasing}%
For all increasing $d,e \in \E$ and $Y = \Ran d \Inter \Ran e$,
$\EqualRes {\Inv d} {\Inv e} Y$.

\ITEM{incr_or_decr}%
For every $e \in \E$, either $e$ is increasing and all the elements in
\Ran e have the same cardinality, or otherwise, $e$ is decreasing and
all the elements in \Dom e have the same cardinality.

\ITEM{minus}%
For all $e \in \E$ and $\xi < \zeta \in \Dom e$, $\zeta - \xi =
e(\zeta) - e(\xi)$.

\end{ITEMS}

\SimpleComment{%
\begin{Proof}%
The properties are obvious consequences of \Definition {E}.
\end{Proof}
}
\end{FACT}

\begin{DEFINITION}{p}
First we need some auxiliary means used in this definition only.  For
all functions $p$ and $e$, $p \FUnion e$ is the function $p \Union
\Par [\big] {\Res e (\Dom e \Minus \Dom p)}$. Let \ReflH be the set of
all limit ordinals below $\kappa$ and \ReflG be the set of all
successor ordinals below $\kappa$. For every $\alpha < \kappa$ let
$\UWo_\alpha$ be a fixed well-ordering of \ULevel \alpha and define a
well-ordering of \U by
 \[
	u \UWo v \Iff\  \Ord u < \Ord v \Or
	(\Ord u = \Ord v = \alpha \And u \UWo_\alpha v).
 \]
For each $w \in \U$, denote the set \Set {w'} {w' \UWo w} by \UPred w.
Furthermore, let \Id u, for $u \in \U$, denote the identity function
 \[
	\Set {\Pair \xi \xi}
	{\xi < \Ord u \And
	\EqualRes {\Fst u} {\Snd u} {\xi+1}}.
 \]

Now define for each $u \in \U$ a function $p_u$ as follows.
\begin{ITEMS}

\ITEM {EmptyTuple}
First of all $p_u = \emptyset$ for $u = \EmptyTuple$.

\ITEM {succ}
Suppose $u \in \Proj \T \TLo$, $u$ is a successor of $v$, and for $v$
the function $p_v$ is already defined. Fix an ordinal $\beta < \Ord u$
as follows:
\begin{enumerate}

\ITEM{succ_H}%
Suppose $u \in \UH$. For all $w \in \U$, define inductively that
\ARRAY[ll]{
	\beta_w' = \min \Par [\Big] {
	\ReflH \Minus & \Par [\big] {
		\Dom {\Id u} \Union
		(\Ord w +1) \\
		&\Union \Set {\beta_w'} {w' \in \UPred w}
	}}.
}%
\PvComment {In this case \Ord u is not a regular cardinal.}
Fix $\beta$ to be $\beta_v'$.

\ITEM {succ_G}%
Suppose $u \in \UG$. By \Definition {G}, $\Ord u \in \GInds$ and there
is a unique $\beta \in \ReflG$ with $\Ord u = \beta + 1$. \Note {Note
that $v = \pi(\Ord u)$.}

\end{enumerate}

Assume \Seq {w^\beta_{\gamma'}} {\gamma' < \theta}, for $\theta \leq
\beth_{\beta+1}$, is a fixed enumeration of \TLevel \beta without
repetition. Let $\gamma$ be the ordinal for which $\Res u \beta =
w^\beta_\gamma$ holds. We define
 \[
    \FunctionDefinition {p_u}{
	\FunctionDefMidCase
		{\Id u \FUnion (p_v \FUnion \EFunc [\gamma] \beta)}
		{\If \Ran {p_v} \Text{is ordinal}\,}
	\FunctionDefOtherwise
		{\Id u \FUnion \Par [\big] {
		p_v \FUnion \Inv [\Par] {\EFunc [\gamma] \beta}}}
    }
 \]

\ITEM {limit}
Suppose $u \in \Proj \T \TLo$, $u$ is a limit, and for all $v \TBelow
u$, functions $p_v$ are defined. Then define $p_u$ to be \BigUnion [v
\TBelow u] {p_v}.  \Note {By the definition of $\FUnion$, $p_w \Subset
p_v$ for all $w \TBelow v \TBelow u$.}

\ITEM {not_TFWo}
For all $u \in \U \Minus \Proj \T \TLo$ define $p_u$ to be \Inv [\Par]
{p_{\Inv u}}.

\end{ITEMS}

For every successor $u \in \U \Minus \Braces [\big] \EmptyTuple$, say
a successor of $v \in \U$, there is unique $e \in \E$ such that either
$e = \emptyset$, or otherwise, $\Dom e = \Dom {p_u} \Minus {\Dom
{p_v}}$ and $p_u = p_v \Union e$. We denote this $e$ by \EFunc u.
\end{DEFINITION}

\Remark The part \Id u in the definition above is needed first time in
\ItemOfLemma {back_and_forth} {step1} to ensure that all the functions
have arbitrary large extensions. Note also that $p_u$ might be $\Id u
= \Set {\Pair \xi \xi} {\xi < \beth_\beta}$ when $\Ord u = \beta +1$,
$u \in \UG$ is a successor of \EmptyTuple, and \EqualRes {\Fst u}
{\Snd u} {\beth_\beta}.

\begin{FACT}{p}
\begin{ITEMS}

\ITEM{id}%
If $u \in \U$, $\xi \in \Dom {p_u}$, and $p_u(\xi) = \xi$, then
$p_u(\zeta) = \zeta$ for all $\zeta \leq \xi$.

\ITEM{dom_ran}%
For every $u \in \U$, $p_u$ is a partial function from $\beth_{\Ord
u}$ into $\beth_{\Ord u}$.

\ITEM{below->subset}%
For all $u,v \in \U$, $u \TBelow v$ implies $p_u \ProperSubset p_v$.

\ITEM{succ}%
For every successor $u \in \U \Minus \Braces \EmptyTuple$, \Dom {p_u}
is the cardinal $\beth_{\Ref {\EFunc u}}$ iff \EFunc u is $\emptyset$
or \EFunc u is increasing, and \Ran {p_u} is the cardinal $\beth_{\Ref
{\EFunc u}}$ iff \EFunc u is $\emptyset$ or \EFunc u is decreasing.

\ITEM{limit}%
For all limit points $u \in \U$, $\Dom {p_u} = \BigUnion [v \TBelow u]
{\Dom {p_v}} = \Ran {p_u} = \BigUnion [v \TBelow u] {\Ran {p_v}} =
\beth_{\Ord u}$. %
\PvComment {When $\Ord u = \beth_{\Ord u}$, $\Dom {p_u} = \Ran {p_u} =
\Ord u$, e.g., for all limit points in \UH.}

\ITEM{cofinality}%
For all limit points $u \in \UG$, \Dom {p_u} is a cardinal of
cofinality less than $\lambda$.

\ITEM{ref}%
Suppose that both $u$ and $v$ are successor elements in \U. If $\EFunc
u \Inter \EFunc v \not= \emptyset$ then $u$ and $v$ are successors of
the same element, $\EFunc u = \EFunc v$, and for $\beta = \Ref {\EFunc
u} = \Ref {\EFunc v}$, \EqualRes u v \beta.

\end{ITEMS}

\begin{Proof}%
The proofs of \Item {dom_ran}--\Item {limit} are straightforward
inductions on the tree order $\TBelow$.  Note that for every limit $u
\in \U$, $u = \BigUnion [v \TBelow u] v$ by \Fact {U}. Note also that
in \ItemOfDefinition {p} {succ_H}, when $u \in \UH$ is a successor of
$v$, the following holds:
 \[
	\beta_v' < \SuccCard [\Par] {\Card {\UPred v}}
	\leq \SuccCard [\Par] {\beth_{\Ord v +1}} < \Ord u,
 \]
since $\Card {\UPred v} \leq \beth_{\Ord v +1}$, $\Ord u = \beth_{\Ord
u}$, and \ReflH is the set of all limit ordinals.

\ProofOfItem{cofinality}%
By \Item {limit} \Dom {p_u} is the cardinal $\beth_{\Ord u}$. By
\ItemOfDefinition {G} {limit}, $\Cf {\Ord u} < \lambda$.

\ProofOfItem{ref}%
Let $w^1$ and $w^2$ be such that $u$ is a successor of $w^1$ and $v$
is a successor of $w^2$. By \ItemOfDefinition {p} {succ}, $p_u =
p_{w^1} \Union \EFunc u$ and $p_v = p_{w^2} \Union \EFunc v$. If $w_1
= w_2$ and $\EFunc u \Inter \EFunc v \not= \emptyset$ then there are
$\beta < \Min {\Ord u, \Ord v}$ and $\gamma < \beth_{\beta +1}$ such
that $\EFunc u = \EFunc v \Subset \EFunc [\gamma] \beta$ and $\Res u
\beta = w^\beta_\gamma = \Res v \beta$, where $w^\beta_\gamma$ is
given in \ItemOfDefinition {p} {succ}.

Suppose $u \in \UH$. If $v \in \UH$ and $w^1 \not= w^2$ then $\Ref
{\EFunc u} \not= \Ref {\EFunc v}$ because the mapping \Mapping w
{\beta_w'}, given in \Definition {G}, is injective. Hence $\EFunc u
\Inter \EFunc v = \emptyset$.  Assume $v \in \UG$. Then $\Ref {\EFunc
v} \in \ReflG \Minus \ReflH$, and because $\Ref {\EFunc u} \in
\ReflH$, $\EFunc u \Inter \EFunc v = \emptyset$ holds. Similarly,
$\EFunc u \Inter \EFunc v = \emptyset$ if $u \in \UG$ and $v \in \UH$.

Suppose both $u \in \UG$ and $v \in \UG$. Then the equations $\pi(\Ord
u) = w^1$ and $\pi(\Ord v) = w^2$ hold, and there are $\beta^1,
\beta^2 \in \ReflG$ with $\Ord u = \beta^1 +1$ and $\Ord v = \beta^2
+1$. If $\Ord u \not= \Ord v$ then $\Ref {\EFunc u} = \beta^1 \not=
\beta^2 = \Ref {\EFunc v}$ and hence $\EFunc u \Inter \EFunc v =
\emptyset$. On the other hand, when \Ord u equals \Ord v, $w^1 =
\pi(\Ord u) = \pi(\Ord v) = w^2$.
\end{Proof}
\end{FACT}
%


%
\begin{DEFINITION}{FComps}
We define \FCompSeqs to be the set of all pairs \Pair {\bar u} {W}
satisfying the following conditions:
\begin{ITEMS}

\ITEM{nonempty}%
$W$ is nonempty.

\ITEM{u}%
For some $n < \omega$, $\bar u$ is a sequence \Seq {u_i} {i < n} of
elements in $\U \Minus \Braces \EmptyTuple$.

\ITEM{composition}%
Let $W_0$ be $W$. Inductively for every $i < n-1$, $W_i \Subset \Dom
{p_{u_i}}$ and $W_{i+1} = \Image {p_{u_i}} {W_i}$.

\ITEM{BegEnd}%
For every $i < n-1$, \EqualRes {\Snd {u_i}} {\Fst {u_{i+1}}} {\PSup
{W_{i+1}}}.

\end{ITEMS}
For every \Pair {\bar u} W in \FCompSeqs there is the natural sequence
\FuncSeq u W of functions defined as follows:
 \[
	\FuncSeq u W =
	\Seq {\FuncSeq [i] u W} {i < \SLh u},
 \]
where each \FuncSeq [i] u W is a shorthand for \Res {p_{u_i}}
{W_i}. The composition $\Func [\SLh u -1] u W \Comp \dots \Comp \Func
[0] u W$ is denoted by \Func u W. For all sequences $\FSeq f = \Seq
{f_i} {i < n}$, $1 \leq n < \omega$, which are of the form \FuncSeq u
W, for some fixed $\Pair {\bar u} W \in \FCompSeqs$ with $\SLh u = n$,
we shall use the following notation:
\begin{itemize}

\item for each $i < \SLh u$,
 \ARRAY{%
	\Ind {f_i} = u_i, \\
	\Beg {f_i} = \Res {\Fst {u_i}} {\PSup {\Dom {f_i}}}, \\
	\End {f_i} = \Res {\Snd {u_i}} {\PSup {\Ran {f_i}}}.
 }

\item $\Beg f = \Beg {f_0}$ and $\End f = \End {f_{\SLh u -1}}$;

\item $f$ is the composition \Composition f 0 {\SLh u -1};

\item for $i < \SLh u$, \ResComp f {\leq i} is a shorthand for
\Composition f 0 i, and for all $\xi \in \Dom f$,
 \[
    \FunctionDefinition {\ResComp f {< i}(\xi)}{
	\FunctionDefMidCase \xi {\If i = 0}
	\FunctionDefOtherwise {\Composition f 0 {i-1}(\xi)}
    }
 \]

\end{itemize}
\end{DEFINITION}

\begin{DEFINITION}{Minimal}%
A sequence $\FSeq f = \FuncSeq u W$, $\Pair {\bar u} W \in
\FCompSeqs$, is called minimal if the following two conditions are
satisfied:
\begin{ITEMS}

\ITEM{dom}%
for all $i < \SLh f$ and $v \TBelow \Ind {f_i}$, $\Dom {f_i}
\not\Subset \Dom {p_v}$;

\ITEM{min}%
there are no indices $i \leq j < \SLh f$ satisfying that the composition
\Composition f i j is identity and $\Beg {f_i} = \End {f_j}$.

\end{ITEMS}
Let \FComps be the set \Set {\FuncSeq u W} {\Pair {\bar u} W \in
\FCompSeqs \And \FuncSeq u W \Text{is minimal\,}}. We abbreviate \Set
{f} {\FSeq f \in \FComps \And \SLh f = 1} to \F.
\end{DEFINITION}

\begin{FACT}{Minimal}%
\begin{ITEMS}

\ITEM{reduction}%
For all $\Pair {\bar u} W \in \FCompSeqs$, either \Func u W is the
identity function and $\Beg {\Func u W} = \End {\Func u W}$, or
otherwise, there is a sequence $\FSeq f \in \FComps$ such that $\SLh f
\leq \SLh u$, $f = \Func u W$, $\Beg f = \Beg {\Func u W}$, and $\End
f = \End {\Func u W}$.

\ITEM{not_id}%
For every $q \in \F$ there is $\theta \in \Dom q$ such that for all
$\xi \in \Dom q \Minus \theta$, $q(\xi) \not= \xi$.

\ITEM{nonempty_efunc}%
For every $q \in \F$ with $\Ind q = u$ a successor, \EFunc u is
nonempty.

\end{ITEMS}

\SimpleComment{%
\begin{Proof}%
We sketch a proof for \Item {reduction}. Suppose \Func u W is not
identity or $\Beg {\Func u W} \not= \End {\Func u W}$.  We have to
show that $\bar u$ can be reduced to a form $\bar {u'}$ such that the
sequence \FuncSeq {u'} W satisfies the conditions \DefinitionItem
{Minimal} {dom} and \DefinitionItem {Minimal} {min} of \Definition
{Minimal}. For each $i < \SLh u$, let $v_i$ be the smallest element
with $v_i \TBelowEq u_i$ and $\Func [i] u W \Subset p_{v_i}$, and
denote \Seq {v_i} {i < \SLh u} by $\bar v$. The pair \Pair {\bar v} W
is in \FCompSeqs since $\Beg {\Func [i] u W} = \Beg {\Func [i] v W}$
and $\End {\Func [i] u W} = \End {\Func [i] v W}$ for every $i < \SLh
u = \SLh v$. The sequence \FuncSeq v W satisfies \ItemOfDefinition
{Minimal} {dom}.

Suppose $i \leq j < \SLh u$ are indices exemplifying that \FuncSeq v W
does not satisfy \ItemOfDefinition {Minimal} {min}. By our assumption
at the beginning of the proof, $i \not= 0$ or $j \not= \SLh v -1$, and
furthermore,
 \[
    \FunctionDefinition {\Beg {\Func [i] v W}}{
	\FunctionDefMidCase   {\End {\Func [i-1] v W}} {\If i > 0}
	\FunctionDefOtherwise [;] {\Beg {\Func u W}}
    }
 \]
and
 \[
    \FunctionDefinition {\End {\Func [j] v W}}{
	\FunctionDefMidCase   {\Beg {\Func [j+1] v W}} {\If j < \SLh v -1}
	\FunctionDefOtherwise {\End {\Func u W}}
    }
 \]
Hence the sequence \FSeq w, which is left over when the elements $v_i,
\dots, v_j$ are removed from \FSeq v, is nonempty. Because of the
equation $\Beg {\Func [i] v W} = \End {\Func [j] v W}$, the pair \Pair
{\bar w} W is in \FCompSeqs, $\Beg {\Func w W} = \Beg {\Func u W}$,
and $\End {\Func w W} = \End {\Func u W}$ in each of the cases $0 = i
\leq j < \SLh v -1$, $0 < i \leq j = \SLh v -1$, and $0 < i \leq j <
\SLh v -1$. So this way all the components of \FSeq v violating
\ItemOfDefinition {Minimal} {min} can be eliminated and the resulting
sequence is as wanted.
\end{Proof}
}
%
\end{FACT}

Note that for all functions $x$ and sets $X$, \Res x X means the
restricted function \Res x {(\Dom x \Inter X)}, so we do not demand
that $X \Subset \Dom x$ in any such restrictions.

\begin{LEMMA}{Res}
For all nonempty $q \in \F$,
\begin{center}%
	\Ind q is a successor if, and only if,
	$\PSup {\Dom q} \not= \PSup {\Ran q}$.
\end{center}
Moreover, if $\Ind q = u$ is a limit, then
 \[
	\PSup {\Dom q} = \PSup {\Ran q}
	= \Dom [\big] {p_u} = \Ran [\big] {p_u}
	= \beth_{\Ord u}.
 \]

\begin{Proof}%
First of all recall, that $q$ is not identity, \ItemOfFact {Minimal}
{not_id}.  Suppose first that $\Ind q = u$ is a successor of $v \in
\U$. Then $q \Subset p_v \Union e$ for $e = \Res {\EFunc u} {\Dom
q}$. Abbreviate \Ref e by $\gamma$. We have $\gamma \geq \Ord v$ and
$\Dom {p_v} \Union \Ran {p_v} \Subset \beth_{\Ord v} \leq
\beth_\gamma$. Because $e$ is nonempty, the claim follows from the
facts that $\Dom e \Inter \beth_\gamma \not= \emptyset$ implies $\Dom
e \Subset \beth_\gamma$ and $\Ran e \Inter \beth_\gamma = \emptyset$,
and on the other hand, $\Dom e \Inter \beth_\gamma = \emptyset$
implies $\Ran e \Subset \beth_\gamma$.

Assume \Ind q is a limit. Denote \Ind q by $u$ and $\beth_{\Ord u}$ by
$\mu$. Because $\Dom {p_u} = \Ran {p_u} = \mu$ it suffices to prove
that $\PSup {\Dom q} \geq \mu$ and $\PSup {\Ran q} \geq \mu$.

Let $\theta < \kappa$ be such that \Seq {v_i} {i < \theta} is a
$\TBelow$-increasing enumeration of the elements $w \TBelow u$. We
know that for all ordinals $i$ in \Set {j + (2n+1)} {j < \theta \Text
{is a limit ordinal or} 0, \And n < \omega}, \Dom {p_{v_i}} is a
cardinal. If $\PSup {\Dom q} < \mu$ then there would be $i < \theta$
with $\Dom q \Subset \Dom {p_{v_i}}$ contrary to \ItemOfDefinition
{Minimal} {dom}. So \Dom q must be unbounded in $\mu$.
\PvComment {\Seq {\Dom {p_{v_i}}} {i < \theta \Text{is odd}} is an
increasing sequence of cardinals with limit $\mu$.}

Besides, we know that for all ordinals $i$ in the set $I = \Set {j +
2n} {j < \theta \Text {is a limit ordinal and} n < \omega}$, \Ran
{p_{v_i}} is a cardinal. So if $\PSup {\Ran q} < \mu$ and $i \in I$ is
such that $\Ran q \Subset \Ran {p_{v_i}}$, then $\Dom q \Subset \Dom
{p_{v_i}}$ since $q$ is injective and \Ran {p_{v_i}} is an
ordinal. Thus \Ran q is also unbounded in $\mu$.
\PvComment{$p_{v_i} \Subset p_u$}
\end{Proof}
\end{LEMMA}

\begin{LEMMA}{FComps}
\begin{ITEMS}

\ITEM{initial}%
Suppose $p, q \in \F$ and $e \in \E$ is a nonempty increasing function
with $e \Subset p \Inter q$ and $\Ref e = \beta$. Then for $X = \Dom p
\Inter \Dom q \Inter \beth_\beta$, \EqualRes p q X and \EqualRes {\Ind
p} {\Ind q} \beta. Particularly, \EqualRes {\Beg p} {\Beg q} {\PSup
X}.

\ITEM{cardinal}%
For all $q \in \F$, if \PSup {\Dom q} is a cardinal, then $\PSup {\Dom
q} \leq \PSup {\Ran q}$.

\ITEM{increasing}%
If $\FSeq f \in \FComps$ and $\PSup {\Dom {f_0}} < \PSup {\Ran {f_0}}$
then $\PSup {\Dom {f_i}} < \PSup {\Ran {f_i}}$ for every $i < \SLh f$.

\ITEM{limit}%
For all $\FSeq f \in \FComps$ and $i < \SLh f$, $\PSup {\Dom {f_i}
\Union \Ran {f_i}} \leq \PSup {\Dom f \Union \Ran f}$.

\ITEM{bounded}%
For all $\FSeq f \in \FComps$, if \PSup {\Ran f} is a cardinal $\mu$
and $\Dom f \Subset \mu$ then $\PSup {\Dom f} = \mu$. Furthermore, if
$\Dom f = \mu$ then $\Dom {f_i} = \Ran {f_i} = \mu$ for each $i < \SLh
f$.

\end{ITEMS}

%
\begin{Proof}%
\ProofOfItem{initial}%
By \Definition p, there are successors $u, v \in \U$ with $u \TBelowEq
\Ind p$, $v \TBelowEq \Ind q$, $\Res {\EFunc u} {\Dom e} = \Res {p_u}
{\Dom e} = e = \Res {p_v} {\Dom e} = \Res {\EFunc v} {\Dom e}$, and
$\Dom {p_u} = \Dom {p_v} = \beth_\beta$. By \ItemOfFact {p} {ref}, $u$
and $v$ are successors of the same element, say $w \in \U$, $\EFunc u
= \EFunc v$, and \EqualRes u v \beta.  Therefore we have
 \[
	\Res p {\beth_\beta} \Subset p_u 
	= p_w \Union \EFunc u = p_w \Union \EFunc v
	= p_v \Superset \Res q {\beth_\beta},
 \]
and
 \[
	\Res {\Beg p} {\beth_\beta}
	\Subset \Res {\Fst u} {\beth_\beta}
	= \Res {\Fst v} {\beth_\beta}
	\Superset \Res {\Beg q} {\beth_\beta}.
 \]

\ProofOfItem{cardinal}%
Suppose that $\PSup {\Dom q} > \PSup {\Ran q}$. Then by \Lemma {Res},
\Ind q must be a successor. Denote \Res {\EFunc {\Ind q}} {\Dom q} by
$d$. Necessarily $d$ is decreasing and \Dom d is an end segment of
\Dom q.  By \ItemOfFact {E} {incr_or_decr}, $\Card [\big] {\PSup {\Dom
d}} = \Card [\big] {\min \Dom d}$. Thus $\PSup {\Dom q} = \PSup {\Dom
d} > \min \Dom d$ is not a cardinal.

\ProofOfItem{increasing}%
Suppose, contrary to the claim, that there is $j \in \Interval 1 {\SLh
f -1}$ with $\PSup {\Dom {f_j}} \geq \PSup {\Ran {f_j}}$. We may
assume that $j$ is the smallest possible index with this property. 

Suppose first that $\PSup {\Dom {f_j}} = \PSup {\Ran {f_j}}$.  Then
for $u = \Ind {f_j}$, \PSup {\Dom {f_j}} is the cardinal $\beth_{\Ord
u}$ by \Lemma {Res}. It follows from the equation $\Dom {f_j} = \Ran
{f_{j-1}} = \Dom {\Inv {f_{j-1}}}$ and by applying \Item {cardinal} to
\Inv {f_{j-1}}, that $\PSup {\Dom {\Inv{f_{j-1}}}} \leq \PSup {\Ran
{\Inv {f_{j-1}}}}$. However, then $\PSup {\Dom {f_{j-1}}} \geq \PSup
{\Ran {f_{j-1}}}$, contrary to the choice of $j$.

So suppose $\PSup {\Dom {f_j}} > \PSup {\Ran {f_j}}$. Note that $\PSup
{\Dom {f_{j-1}}} < \PSup {\Ran {f_{j-1}}}$. Abbreviate \Ind {f_{j-1}}
by $u$ and \Ind {f_j} by $v$. By \Lemma {Res}, there are $w^1, w^2 \in
\U$ such that $u$ is a successor of $w^1$ and $v$ is a successor of
$w^2$. Denote $\Ref {\EFunc u}$ by $\gamma^1$, $\Ref {\EFunc v}$ by
$\gamma^2$, \Res {\EFunc u} {\Dom {f_{j-1}}} by $d^1$, and \Res
{\EFunc v} {\Dom {f_j}} by $d^2$. Then $d^1$, $d^2$ are nonempty,
$d^1$ is increasing, $d^2$ is decreasing, and
 \ARRAY{
	f_{j-1} \Subset p_u = p_{w^1} \Union d^1, \\
	\Ran{p_{w^1}} \Subset \beth_{\gamma^1}, \\
	\Ran {d^1} \Subset \beth_{\gamma^1 +1} \Minus
		\beth_{\gamma^1}, \\
	f_j \Subset p_v = p_{w^2} \Union d^2 \\
	\Dom {p_{w^2}} \Subset \beth_{\gamma^2}, \\
	\Dom {d^2} \Subset \beth_{\gamma^2 +1} \Minus
		\beth_{\gamma^2}. \\
 }
Because $\Ran {f_{j-1}} = \Dom {f_j}$, it follows that $\gamma_1 =
\gamma_2$ and $\Ran {d^1} = \Dom {d^2}$. By \ItemOfFact {E}
{increasing}, $d^1 = \Inv {(d^2)}$, and hence $\EFunc u \Inter \EFunc
{(\Inv v)} \not= \emptyset$ \Note {the notation \Inv v is explained in
\Definition {TLo}}. By \ItemOfFact {p} {ref}, $w^1 = \Inv {(w^2)}$ and
\EqualRes u {\Inv v} {\gamma_1}. Consequently, $f_{j-1} = \Inv {f_j}$
and for $\theta = \PSup {\Dom {f_{j-1}}} = \PSup {\Ran {f_j}} \leq
\beth_{\gamma_1}$, $\Beg {f_{j-1}} = \Res {\Fst u} \theta = \Res {\Fst
{\Inv v}} \theta = \Res {\Snd v} {\theta} = \End {f_j}$ contrary to
the minimality of \FSeq f.

\ProofOfItem{limit}%
Denote \PSup {\Dom f \Union \Ran f} by $\theta$.  To reach a
contradiction let $j < \SLh f$ be the smallest index with $\PSup {\Dom
{f_j} \Union \Ran {f_j}} > \theta$.  Then $\PSup {\Dom {f_j}} \leq
\theta < \PSup {\Ran {f_j}}$ since $\Dom {f_0} = \Dom f$ and $\Ran
{f_i} = \Dom {f_{i+1}}$ for every $i < j$. However, from \Item
{increasing} it follows that $\PSup {\Ran {f_j}} \leq \PSup {\Dom
{f_k}} < \PSup {\Ran {f_k}}$ for all $k \in \Interval {j+1} {\SLh f
-1}$, and so $\PSup {\Ran f} = \PSup {\Ran {f_{\SLh f -1}}} > \theta$,
a contradiction.

\ProofOfItem{bounded}%
By applying \Item {cardinal} to \Inv {f_n}, for $n = \SLh f -1$, we
get that $\PSup {\Dom {f_n}} \geq \PSup {\Ran {f_n}} = \mu$. By \Item
{limit}, $\PSup {\Dom {f_n}} \leq \mu$. Hence $\PSup {\Dom {f_n}} =
\mu$. In the same way it can be shown $\PSup {\Dom {f_i}} = \mu$ for
all $i \leq n$.

Suppose $\Dom f = \mu$. By \Lemma {Res}, \Ind {f_i} is a limit point,
say $u_i \in \U$, and $\PSup {\Ran {f_i}} = \mu = \Dom {p_{u_i}} =
\Ran {p_{u_i}}$ for every $i \leq n$. However, when $n > 1$, $f_0
\Subset p_{u_0}$ together with $\Dom {f_0} = \Dom f = \mu = \Dom
{p_u}$ imply that $f_0 = p_{u_0}$ and $\Ran {f_0} = \mu = \Dom
{f_1}$. A similar reasoning shows $\Dom {f_i} = \Ran {f_i} = \mu$ for
every $i \leq n$.
\end{Proof}%
%
%
\end{LEMMA}

\begin{LEMMA}{Ran}%
\begin{ITEMS}

\ITEM{unique}%
Suppose $\FSeq f, \FSeq g \in \FComps$ are such that $\Dom f = \Dom g
= \Braces \xi$, $f(\xi) = g(\xi)$, and $f_i$ is increasing for every
$i < \SLh f$. Then $\SLh f \leq \SLh g$ and for $k = \SLh g - \SLh f$,
both $f_i = g_{k + i}$ and $\Beg {f_i} = \Beg {g_{k + i}}$ hold for
every $i < \SLh f$.  Moreover, if $\Beg g = \Beg f$ or for every $j <
\SLh g$, $g_j$ is increasing, then $\SLh f = \SLh g$.

\ITEM{card}%
Suppose $q \in \F$, $\xi < \zeta \in \Dom q$, and that both \Res q
{\Braces \xi} and \Res q {\Braces \zeta} are increasing. If there is
no $d \in \E$ with $\Braces {\xi, \zeta} \Subset \Dom d$ then $\Card
\xi < \beth_{\Ref {\Res q {\Braces \xi}}} \leq \Card \zeta$.

\ITEM{step}%
Suppose $\SimpleSeq {q_1, q_2} \in \FComps$ and $\PSup {\Dom {q_1}} <
\PSup {\Ran {q_1}}$. Then \Ind {q_1}, \Ind {q_2} are successors, the
functions $d^1 = \Res {\EFunc {\Ind {q_1}}} {\Dom {q_1}}$, $d^2 = \Res
{\EFunc {\Ind {q_2}}} {\Dom {q_2}}$ satisfy the demand that they both
are increasing, and $\Dom {q_2} \Minus \min \Ran {d^1} \Subset \Dom
{d^2}$.

\ITEM{efunc}%
Suppose $\FSeq f \in \FComps$ and $\PSup {\Dom {f_0}} < \PSup {\Ran
{f_0}}$. Then for every $i < \SLh f$, $\Ind {f_i}$ is successor,
\EFunc {\Ind {f_i}} is increasing, and particularly, for $d = \Res
{\EFunc {\Ind {f_0}}} {\Dom {f_0}}$ and for every $i \in \Interval 1
{\SLh f -1}$, $\Dom {f_i} \Minus \ResComp f {<i} \Par [\big] {\min
\Ran d} \Subset \EFunc {\Ind {f_i}}$.

\ITEM{seq}%
Let \FSeq f be as in \Item {efunc}. For $\bar u = \Seq {\Ind {f_i}} {i
< \SLh f}$ and for every $\theta$ with $\PSup {\Dom f} < \theta \leq
\beth_{\Ref d}$, the pair \Pair {\bar u} \theta is in \FCompSeqs and
the sequence \FuncSeq u \theta is in \FComps \Note {see \Definition
{FComps}}. Furthermore, if \End f is a function with a constant value,
then $\End {\Func u \theta} \Superset \End f$ is also a constant
function.

\ITEM{limit}%
For every \FSeq f in \FComps there is $\xi \in \Dom f$ with $\Ran f
\Subset f(\xi) + \PSup {\Dom f}$.

\ITEM{ran}%
Suppose $\FSeq g \in \FComps$, $\Dom g$ is a cardinal $\mu$, $\Ind
{g_0} = u_0$ is successor, and \EFunc {u_0} is increasing. Assume
$\FSeq h \in \FComps$ is such that $\Dom h = \Braces {\xi, \xi'}
\Subset \mu$, $\xi \in \Dom {\EFunc {u_0}}$, $h(\xi) = g(\xi)$, and
$\Beg h \Subset \Beg g$. Then either $h(\xi') \in \Ran g$, or
otherwise, $h(\xi') \geq h(\xi) + \mu$.

\end{ITEMS}

%
\begin{Proof}%
\ProofOfItem{unique}%
Denote $\SLh f -1$ by $n$ and $\SLh g -1$ by $m$.  If $g_m$ was
decreasing, then, by applying \ItemOfLemma {FComps} {increasing} to
the sequence \InvFSequence g 0 m in \FComps, $g_i$ should be
decreasing for every $i \leq m$ and $\PSup {\Dom g} > \PSup {\Ran
{g_m}} = \PSup {\Ran {f_n}} > \PSup {\Dom f} = \xi +1$ contrary to the
assumption $\Dom f = \Dom g = \Braces {\xi}$. Thus $g_m$ is
increasing. Since $\Ran {f_n} = \Ran {g_m}$, $f_n = g_m$ by
\ItemOfFact {E} {increasing}. By \ItemOfLemma {FComps} {initial},
$\Beg {f_n} = \Beg {g_m}$. When $n > 0$, $\Ran {f_{n-1}} = \Dom {f_n}
= \Dom {g_m} = \Ran {g_{m-1}}$.  Hence we can repeat the same argument
and we get that $f_{n-i} = g_{m-i}$ and $\Beg {f_{n-i}} = \Beg
{g_{m-i}}$ for every $i \leq \Min {m, n}$. However $m \geq n$ since
otherwise $\Dom g = \Braces {\ResComp f {<n-m} (\xi)} \not= \Braces
\xi$.

If $m > n$ and $\Beg g = \Beg f$ then $\ResComp g {\leq m-n-1} (\xi) =
\xi$ and $\End {g_{m-n-1}} = \Beg {g_{m-n}} = \Beg {f_0} = \Beg f =
\Beg g = \Beg {g_0}$ contrary to the minimality of \FSeq g.

If $m > n$ and $g_j$ is increasing, for every $j < \SLh g$, then $\Ran
{g_{m-n-1}} = \Dom {g_{m-n}} = \Dom {f_0} = \Braces \xi$ and $\Dom g =
\Braces {\InvComposition g 0 {m-n-1}(\xi)} \not= \Braces {\xi}$, a
contradiction.

\ProofOfItem{card}%
Let $u \TBelowEq \Ind q$ be the smallest element with $\xi \in \Dom
{p_u}$, and $v \TBelowEq \Ind q$ be the smallest element with $\zeta
\in \Dom {p_v}$. Then $u$ and $v$ are successors, $u \TBelowEq v$,
$\xi \in \Dom {\EFunc u}$, $\zeta \in \Dom {\EFunc v}$, and $\Res q
{\Braces {\xi, \zeta}} \Subset \EFunc u \Union \EFunc v \Subset p_u
\Union \EFunc v$.  Assume $\EFunc u \not= \EFunc v$. Then $u \TBelow
v$.  Since \EFunc u is increasing $\Dom {p_u} = \beth_\beta$ where
$\beta = \Ref {\EFunc u}$. Because $\zeta \in \Dom {\EFunc v} \Minus
\Dom {\EFunc u}$ we have $\Dom {\EFunc v} \Inter \beth_\beta =
\emptyset$. So $\Card \xi < \beth_\beta \leq \Card \zeta$.

\ProofOfItem{step}
follows that $\PSup {\Dom {q_2}} < \PSup {\Ran {q_2}}$. The elements
$\Ind {q_1}$ and $\Ind {q_2}$ are successors by \Lemma {Res}. If
$q_2(\xi) = \xi$ then $q_2$ should be identity contrary to \ItemOfFact
{Minimal} {not_id}. So both $d^1$ and $d^2$ are increasing. If for
some $\xi \in \Ran {d^1}$, \Res {q_2} {\Braces \xi} is decreasing,
then $\Res {\Inv {(d^1)}} {\Braces \xi} = \Res {q_2} {\Braces \xi}$,
and as in the proof of \ItemOfLemma {FComps} {increasing}, $\Inv
{(q_1)} = q_2$ and $\Beg {q_1} = \End {q_2}$ contrary to the
minimality of the sequence \SimpleSeq {q_1, q_2}. Thus \Res {q_2}
{\Braces \xi} is increasing for all $\xi \in \Ran {d^1}$.  Now $\Card
\xi = \Card {\min \Ran {d^1}}$ for each $\xi \in \Ran {d^1}$.  By
\Item {card}, there is $e \in \E$, with $e \Subset q_2$ and $\Dom e =
\Ran {d^1}$. Because $\Ran {d^1} = \Dom e$ is an end segment of $\Ran
{q_1} = \Dom {q_2}$, $\Dom {q_2} \Minus \min \Ran {d^1} = \Dom
e$. Since $e$ is increasing, it follows from the definition of $d^2$
that $e \Subset d^2$.

\ProofOfItem{efunc}%
The claim follows from \Item {step} by induction on $i < \SLh
f$.

\ProofOfItem{seq}%
It suffices to show that $\Pair {\bar u} \theta \in \FCompSeqs$ since
then the minimality of \FuncSeq u \theta follows from the fact that
\FSeq f is in \FComps. Let $\beta_i$ denote \Ref {\EFunc {u_i}} for
every $i < \SLh u = \SLh f$. Abbreviate $\min \Dom d$ by $\xi$. By
\Item {efunc}, $\Ran {p_{u_i}} \Subset \beth_{(\beta_i) +1} \leq
\beth_{\beta_{(i+1)}} = \Dom {p_{u_{i+1}}}$ and $\ResComp f {\leq i}
(\xi) \in \beth_{\beta_i +1} \Minus \beth_{\beta_i}$ for every $i <
\SLh u -1$. So \Pair {\bar u} \theta satisfies \ItemOfDefinition
{FComps} {composition}. From $\Pair {\bar u} {\Dom f} \in \FCompSeqs$
it follows that \EqualRes {\Snd {u_i}} {\Fst {u_{i+1}}} {\Par [\big]
{\ResComp f {\leq i} (\xi) +1}} for all $i < \SLh u -1$. These
equations together with \Definition {U} ensure that both of the
functions \Snd {u_i} and \Fst {u_{i+1}}, for $i < \SLh u -1$, have the
same constant value on the interval $\beth_{\beta_i +1} \Minus
\beth_{\beta_i}$. Hence the pair \Pair {\bar u} \theta satisfies
\ItemOfDefinition {FComps} {BegEnd}, too. Similarly, the latter claim,
concerning \End f, is a consequence of the facts that for $n = \SLh u
-1$, $f(\xi) \in \Ran {\Func u \theta} = \Ran {\Func [n] u \theta}
\Subset \beth_{\beta_n +1} \Minus \beth_{\beta_n}$ and \Snd {u_n} is a
constant function on the interval $\beth_{\beta_n +1} \Minus
\beth_{\beta_n}$.

\ProofOfItem{limit}%
Abbreviate \PSup {\Dom f} by $\theta$, $\SLh f - 1$ by $n$, and for
every $i \leq n$, \Ind {f_i} by $u_i$. If $\PSup {\Ran f} \leq \theta$
there is nothing to prove. So assume $\PSup {\Ran f} > \theta$.  There
must be the smallest index $j \leq n$ satisfying $\PSup {\Dom {f_j}}
\leq \theta < \PSup {\Ran {f_j}}$, $u_j$ is a successor, and \Res
{\EFunc {u_j}} {\Dom {f_j}}, abbreviated by $d$, is increasing. Let
$\xi$ be $\min \Dom d$. Then for all $\zeta \in \Dom {f_j} \Minus \Dom
d$, $f_j(\zeta) < f_j(\xi)$, by the definition of \EFunc {u_j}. %
\PvComment{\Res {f_j} {\Braces \zeta} must be decreasing if $\zeta >
\xi$.}
Besides, $\Dom {f_j} \Subset \theta$ together with \ItemOfFact {E}
{minus} ensure that for all $\zeta \in \Dom d$, $f_j(\zeta) - f_j(\xi)
= d(\zeta) - d(\xi) = \zeta - \xi < \theta$. So the claim holds in
case $j = n$.

Suppose $n > j$. From \Item {efunc} it follows that for every $i \in
\Interval {j+1} n$, $u_i$ is a successor, \EFunc {u_i} is increasing,
and $\Dom {f_i} \Minus \ResComp f {<i} (\xi) \Subset \Dom {\EFunc
{u_i}}$. For all $\zeta \in \Dom {f_j} \Minus \Dom d$, $f(\zeta) <
f(\xi)$ since $f_j(\zeta) < f_j(\xi)$ and the property ``\Res {f_i}
{\Braces \xi} is increasing for every $i \in \Interval {j+1} n$''
implies $\ResComp f {\leq i} (\zeta) < \ResComp f {\leq i} (\xi)$ for
every $i \in \Interval {j+1} n$. Suppose $\zeta \in \Dom d$, $i \in
\Interval {j+1} n$, and $\ResComp f {< i} (\xi) < \ResComp f {< i}
(\zeta) < \ResComp f {< i} (\xi) + \theta$. Then $\Braces {\ResComp f
{<i} (\xi), \ResComp f {< i} (\zeta)} \Subset \Dom {\EFunc {u_i}}$ and
by \ItemOfFact {E} {minus},
 \[
	\ResComp f {\leq i} (\zeta) - \ResComp f {\leq i} (\xi) =
	\EFunc {u_i} (\ResComp f {< i} (\zeta)) -
		\EFunc {u_i} (\ResComp f {< i} (\xi)) =
	\ResComp f {< i} (\zeta) - \ResComp f {< i} (\xi)
	< \theta.
 \]
The claim follows from the fact that $\Ran f \Minus f(\xi) = \Ran
{f_n} \Minus \ResComp f {\leq n} (\xi) \Subset \Ran {\EFunc {u_n}}$
\Note {remember $\Dom {f_n} \Minus \ResComp f {<n} (\xi) \Subset \Dom
{\EFunc {u_n}}$ and \EFunc {u_n} is increasing}.

\ProofOfItem{ran}%
Denote $\SLh g -1$ by $n$, $\SLh h -1$ by $m$ and for each $i \leq m$
abbreviate $\ResComp h {<i} (\xi)$ by $\xi_i$ and $\ResComp h {<i}
(\xi')$ by $\xi_i'$. Write $\xi_{m+1}$ for $h(\xi)$ and $\xi_{m+1}'$
for $h(\xi')$. Note that by \Item {efunc}, for every $i \leq n$ and
for $u_i = \Ind {g_i}$, $u_i$ is a successor, \EFunc {u_i} is
increasing, and \EqualRes {g_i} {\EFunc {u_i}} {\Braces {\ResComp g
{<i} (\xi)}}.

There exists the smallest index $j \leq m$ with $\Res {h_j} {\Braces
{\xi_j}} = \Res {g_0} {\Braces {\xi}}$, because otherwise for the
minimal reduct \FSeq d of the sequence \Seq {\Res {h_i} {\Braces
{\xi_i}}} {i \leq m} \Note {see \Fact {Minimal}}, \FSeq d differs from
the minimal sequence $\FSeq e = \Seq {\Res {g_i} {\Braces {\ResComp g
{<i} (\xi)}}} {i \leq n}$ and $\Composition e 0 n = \Composition d 0
{\SLh d -1}$ contrary to \Item {unique}. We have two cases to
consider:
\begin{enumerate}

\ITEM {greater} $\xi_j' \geq \mu$;

\ITEM {smaller} $\xi_j' < \mu$.

\end{enumerate}

\ProofOfItem {greater}%
Suppose first that $\xi_j' \geq \mu$. Note that $\mu > \xi =
\xi_j$. Note that $h_j(\xi_j') \not= \xi_j'$ since otherwise also
$h_j(\xi_j) = \xi_j$.  The function \Res {h_j} {\Braces {\xi_j'}} must
be increasing, namely otherwise, we reach a contradiction in the
following manner.  Assume \Res {h_j} {\Braces {\xi_j'}} is
decreasing. There are two subcases:

\begin{itemize}

\item[i)]%
Assume that $\PSup {\Dom {h_j}} > \PSup {\Ran {h_j}}$ or $\PSup {\Dom
{h_{j-1}}} > \PSup {\Ran {h_{j-1}}}$.

\item[ii)]%
By \Lemma {Res}, both $\PSup {\Dom {h_j}} \not= \PSup {\Ran {h_j}}$
and $\PSup {\Dom {h_{j-1}}} \not= \PSup {\Ran {h_{j-1}}}$ hold. So
suppose $\PSup {\Dom {h_j}} < \PSup {\Ran {h_j}}$ and $\PSup {\Dom
{h_{j-1}}} < \PSup {\Ran {h_{j-1}}}$. 

\end{itemize}

i) It would follow from the assumption $\xi_j' \geq \mu$ and by
applying \Item {efunc} to the sequence \InvFSequence h 0 j or
\InvFSequence h 0 {j-1}, that $\PSup {\Dom h} = \PSup {\Dom {h_0}} >
\mu$, a contradiction.

ii) The function \Res {h_{j-1}} {\Braces {\xi_{j-1}'}} is increasing,
since otherwise,
 \[
	\PSup {\Dom {h_{j-1}}} \geq \xi_{j-1}' +1 \geq h_{j-1}
	(\xi_{j-1}') +1 = \xi_j' +1 = \PSup {\Ran {h_{j-1}}}.
 \]
Let $\beta$ be \Ref {\Res {h_{j-1}} {\Braces {\xi_{j-1}'}}}.  If
$\xi_{j-1} \geq \beth_\beta$, then $\xi_{j-1} > \xi_{j-1}'$ and
$h(\xi_{j-1}) \not= \xi_{j-1}$. Moreover, $\Ref {\Res {h_{j-1}}
{\Braces {\xi_{j-1}}}} > \beta$ and $\PSup {\Dom {h_{j-1}}} =
\xi_{j-1} +1 > \beth_{\beta +1} > \xi_j' +1 = \PSup {\Ran {h_{j-1}}}$,
a contradiction. On the other hand, if $\xi_{j-1} < \beth_\beta$, then
it follows from the assumption $\xi_j' > h_j (\xi_j')$ that $\Res
{h_{j-1}} {\Braces {\xi_{j-1}'}} = \Inv {\Par [\big] {\Res {h_j}
{\Braces {\xi_j'}}}}$. By \ItemOfLemma {FComps} {initial}, $h_j = \Inv
{h_{j-1}}$ and $\Beg {h_{j-1}} = \End {h_j}$ contrary to the
minimality of \FSeq h.

Hence \Res {h_j} {\Braces {\xi_j'}} is increasing, and by \Item
{efunc}, \Ind {h_i}, abbreviated by $v_i$, is successor, \EFunc {v_i}
is increasing, and $\xi_i' \in \Dom {\EFunc {v_i}}$ for every $i \in
\Interval j m$. We show by induction on $i \in \Interval j m$ that
$\xi_i + \mu \leq \xi_i'$ where $+$ is the ordinal addition.

Since $\xi_j = \xi < \mu$, $\xi_j' \geq \mu$, and $\mu$ is cardinal,
we have $\xi_j + \mu \leq \xi_j'$. Suppose $i < m$ and $\xi_i + \mu
\leq \xi_i'$. If $h_i(\xi_i) = \xi_i = \xi_{i+1}$ then $\xi_{i+1} +
\mu = \xi_i+\mu \leq \xi_i' < h_i(\xi_i') = \xi_{i+1}'$.

If $\xi_i \in \Dom {\EFunc {v_i}}$, then $\xi_{i+1}' = \EFunc
{v_i}(\xi_i')$, $\xi_{i+1} = \EFunc {v_i}(\xi_i)$, and $\xi_{i+1}' -
\xi_{i+1} = \EFunc {v_i}(\xi_i') - \EFunc {v_i}(\xi_i) = \xi_i' -
\xi_i \geq \mu$.

If $\xi_i \not\in \Dom {\EFunc {v_i}}$, then the reflection point of
\Res {h_i} {\Braces {\xi_i}}, say $\beta$, is smaller than \Ref
{\EFunc {v_i}} by the definition of \EFunc {v_i}.  Since $\xi_{i+1}
\leq \beth_{\beta+1} \leq \beth_{\Ref {\EFunc {v_i}}}$ and $\mu <
\beth_{\Ref {\EFunc {v_i}}}$, it follows that $\xi_{i+1} + \mu <
\beth_{\Ref {\EFunc {v_i}}} < \xi_{i+1}'$.

\ProofOfItem {smaller}%
Suppose then that $\xi_j' < \mu$.  Abbreviate \Ref {\EFunc {u_i}}, for
$i \leq n$, by $\gamma_i$.  Since $\Res {h_j} {\Braces {\xi_j}} = \Res
{\EFunc {u_0}} {\Braces {\xi}} = \Res {g_0} {\Braces \xi}$ is
increasing and $\Braces {\xi_j, \xi_j'} \Subset \mu = \Dom {g_0} \leq
\beth_{\gamma_0}$, we get by \ItemOfLemma {FComps} {initial}, that
$\Beg {h_j} \Subset \Beg {g_0}$ and \EqualRes {h_j} {g_0} {\Braces
{\xi_j'}}. By \Item {efunc}, \Res {h_i} {\Braces {\xi_i}} are
increasing for all $i \in \Interval j m$. It follows from $\xi =
\xi_j$ and $h(\xi_j) = h(\xi) = g(\xi)$ together with
\Item {unique}, that $\Seq {\Res {h_i} {\Braces {\xi_i}}} {j \leq i
\leq m} = \Seq {\Res {g_k} {\Braces {\ResComp g {<k} (\xi)}}} {k \leq
n}$.

To show that $h(\xi') = \Composition h j m (\xi_j') \in \Ran g$ we
prove by induction on $k \leq n$ that \EqualRes {h_{j+k}} {g_k}
{\Braces {\xi_{j+k}'}}. Note that $m = j+n$ and it is possible that
$\xi_j' \not= \xi'$. We already proved the case $k = 0$. Suppose $k >
0$ and for every $i < k$ the subclaim holds. Then $\Braces {\xi_{j+k},
\xi_{j+k}'} = \Dom {h_{j+k}} \Subset \Ran {g_{k-1}} = \Dom
{g_k}$. Since $\Res {h_{j+k}} {\Braces {\xi_{j+k}}} = \Res {g_k}
{\Braces {\ResComp g {<k} (\xi)}}$ is increasing and $\xi_{j+k}' \in
\Dom {g_k} \Subset \beth_{\gamma_k}$, we get by \ItemOfLemma {FComps}
{initial}, that \EqualRes {h_{j+k}} {g_k} {\Braces {\xi_{j+k}'}}.
\end{Proof}%
%
%
\end{LEMMA}

\begin{LEMMA}{Limit}%
\begin{ITEMS}%
\ITEM{increasing}%
Suppose $p,q \in \F$ are such that \Ind p is a limit, $\Dom p = \Dom q
= X$, and the set $Y = \Set {\zeta \in X} {\Res p {\Braces \zeta} =
\Res q {\Braces \zeta}}$ is unbounded in $X$. Then $\Ind q = \Ind p$,
and particularly, $p = q$, $\Beg p = \Beg q$, and $\End p = \End q$.

\ITEM{end_seg}%
Suppose \FSeq f in \FComps and the set $I_0 = \Set {\xi \in \Dom
{f_0}} {\xi < f_0(\xi)}$ is unbounded in \Dom {f_0}. Then there is an
end segment $J$ of $I_0$ such that for every $\xi \in J$ and $i < \SLh
f$, \Res {f_i} {\Braces {\ResComp f {<i} (\xi)}} is increasing.

\ITEM{componentwise}%
Suppose $\FSeq f, \FSeq g \in \FComps$ and $n < \omega$ are such that
$\SLh f = \SLh g = n$, $\Dom f = \Ran f$ is a cardinal, and $f \Subset
g$. Then $\Ind {f_i} \TBelowEq \Ind {g_i}$ and $f_i \Subset g_i$ holds
for every $i < n$.

\ITEM{id}%
Suppose $\bar \alpha$ is an increasing sequence \Seq {\alpha_l} {l <
\omega} of ordinals below $\kappa$ such that \PSup {\bar \alpha} is a
cardinal. Then for every $\FSeq f \in \FComps$ such that $\Dom f =
\Set {\alpha_l} {l < \omega}$, there are infinitely many indices $l <
\omega$ with $f(\alpha_l) \not= \alpha_l$.

\end{ITEMS}

%
\begin{Proof}%
\ProofOfItem{increasing}%
Let $u$ be \Ind p and $v$ be \Ind q.  We may assume that $p(\zeta) =
q(\zeta) \not= \zeta$ for all $\zeta \in Y$.  Let $Z$ be the set \Set
{\Min {\zeta, p(\zeta)}} {\zeta \in Y}.  By \ItemOfLemma {FComps}
{initial}, \EqualRes u v {\xi} for every $\xi \in Z$. Since $Y$ is
unbounded in $X$ and \Ind p is a limit, also $Z$ is unbounded in $X$.
By \Lemma {Res}, $\PSup X = \PSup {\Dom p} = \Dom {p_u} = \beth_{\Ord
u}$. Hence, as in the proof of \Fact {U}, $u \TBelowEq v$. By
\ItemOfDefinition {Minimal} {dom}, $u = v$. Because $p$ and $q$ have
the common domain $X$, it follows that $p = q$, $\Beg p = \Beg q$ and
$\End p =
\End q$.

\ProofOfItem{end_seg}%
If \Ind {f_0} is a successor then $\PSup {\Dom {f_0}} < \PSup {\Ran
{f_0}}$, and the claim follows from \ItemOfLemma {Ran} {efunc}.
Suppose \Ind {f_0} is a limit.  By \Lemma {Res}, \PSup {\Dom {f_0}} is
a cardinal. Suppose, contrary to the claim, that $j < n$ is the
smallest index for which there is unbounded $J \Subset I_0$ such that
for every $\zeta$ in the set $Y = \Set {\ResComp f {<j} (\xi)} {\xi
\in J}$, \Res {f_j} {\Braces \zeta} is increasing and \Res {f_{j+1}}
{\Braces {f_j(\zeta)}} is not increasing. Then \Ind {f_i} is a limit
for every $i \leq j$, since otherwise the existence of the chosen $j$
contradicts \ItemOfLemma {Ran} {efunc}. Therefore $\PSup {\Dom {f_j}}$
is the cardinal \PSup {\Dom {f_0}}, and necessarily $Y$ is unbounded
in \Dom {f_j}. So we may assume that \Res {f_{j+1}} {\Braces
{f_j(\zeta)}} is decreasing for every $\zeta \in Y$. Since \Res {f_j}
{\Braces \zeta} must equal \Inv{\Par [\big] {\Res {f_{j+1}} {\Braces
{f_j(\zeta)}}}} for every $\zeta \in Y$ it follows from \Item
{increasing} that $f_j = \Inv {f_{j+1}}$ and $\Beg {f_j} = \End
{f_{j+1}}$ contrary to the minimality of \FSeq f.

\ProofOfItem{componentwise}%
In the case $n = 1$ the claim is proved in \Item {increasing}. Assume
$n > 1$.  Let $\theta$ be the cardinal $\Dom f = \Ran f$.  For each $i
< n$, $\Dom {f_i} = \Ran {f_i} = \theta$ by \ItemOfLemma {FComps}
{bounded}. By \Lemma {Res}, \Ind {f_i} is a limit point and $f_i =
p_{\Ind {f_i}}$ for every $i < n$. Denote the set \Set {\xi < \theta}
{f_0 (\xi) > \xi} by $I_0$. Then $I_0$ must be unbounded in $\theta$
by \Definition {p}. For each $i < n - 1$ define $I_{i +1}$ to be \Set
{\xi \in \Image {f_i} {I_i}} {f_{i +1} (\xi) > \xi}.

By \Item {end_seg}, there is an end segment $K$ of $I_0$ satisfying
that for every $i < n -1$, \Image {\ResComp f {\leq i}} K is an end
segment of $I_{i+1}$. Now $\SLh f = \SLh g = n$ and $f(\xi) = g(\xi)$
together with \ItemOfLemma {Ran} {unique} imply that for all $\xi \in
K$ and $i < n$, the equations $\ResComp f {\leq i}(\xi) = \ResComp g
{\leq i}(\xi)$ hold. Since $K$ is unbounded in $\theta$ and for each
$i < n$, \Res {f_i} {\Image {\ResComp f {<i}} K} is increasing, also
\Image {\ResComp f {\leq i}} K is unbounded in $\theta$ for every $i <
n$. By \Item {increasing}, $\Ind {f_i} \TBelowEq \Ind {g_i}$ and $f_i
= \Res {f_i} \theta = \Res {g_i} \theta$ for every $i < n$.

\ProofOfItem{id}%
Let $\theta$ be the cardinal \PSup {\bar \alpha} and let $n$ denote
the length of \FSeq f. For every $l < \omega$ and $i < n$ write
$d^l_i$ for the function \Res {f_i} {\Braces [\big] {\ResComp f {<i}
(\alpha_l)}}. For every $i < n$ define $I_i$ to be the set \Set {l <
\omega} {d^l_i \Text{is increasing}}.

Suppose, contrary to the claim, that there is $m < \omega$ such that
$f(\alpha_l) = \alpha_l$ hold for all $l \in \omega \Minus m$.  By
\Item {end_seg}, $I_0$ is finite. There must be the smallest $j \in
\Interval 1 {n-1}$ such that $I_j$ is infinite. By \ItemOfLemma
{FComps} {cardinal}, $\PSup {\Dom {f_0}} = \PSup {\Ran {f_0}}$ and so
\Ind {f_0} is a limit.  Since $\Dom {f_{i+1}} = \Ran {f_i}$ for all $i
< j$, we get by applying \ItemOfLemma {FComps} {cardinal} repeatedly,
that \Ind {f_i} are limit points for all $i < j$. By the choice of
$j$, there is an end segment $J$ of $\omega$ such that $\min J \geq m$
and $d^l_i$ is decreasing for all $l \in J$ and $i < j$ \Note {$d^l_i$
cannot be identity for unbounded many $l < \omega$}. The set $Y = \Set
{\ResComp f {<j} (\alpha_l)} {l \in J}$ is unbounded in \Dom {f_j}.

If \Ind {f_j} is a successor, then $Y \Inter \Dom {\EFunc {\Ind
{f_j}}}$ is infinite, and by \ItemOfLemma {Ran} {efunc}, $f(\alpha_l)
> \theta > \alpha_l$ for infinitely many $l \in J$, a
contradiction. Hence \Ind {f_j} is a limit. By \Item {end_seg}, there
is an end segment $K$ of $J$ such that for every $l \in J$ and $k \in
\Interval j {n-1}$, $d^l_k$ is increasing. For every $l \in K$,
$f(\alpha_l) = \alpha_l$ holds, and thus the compositions
$\Inv{(d^l_0)} \Comp \dots \Comp \Inv{(d^l_{j-1})}$ and $d^l_{n-1}
\Comp \dots \Comp d^l_j$ are equal.  Since $\End {f_{j-1}} = \Beg
{f_j}$ and the sequences \SimpleSeq{\Inv{(d^l_{j-1})}, \dots,
\Inv{(d^l_0)}} and \SimpleSeq {d^l_j, \dots, d^l_{n-1}} are in \FComps
\Note {in both of the sequences all the functions are increasing}, it
follows from \ItemOfLemma {Ran} {unique}, that these sequences are
equal. Particularly, $d^l_j = \Inv{(d^l_{j-1})}$ for every $l \in
K$. From \Item {increasing} it would follow that $f_{j-1} = f_j$ and
$\Beg {f_{j-1}} = \End {f_j}$, contrary to the minimality of \FSeq f.
\end{Proof}%
%
%
\end{LEMMA}
%


\begin{DEFINITION}{D}
We define \D to be the following closed and unbounded subset of
$\kappa$:
 \[
	\Set [\Big] {\mu \in \kappa \Minus (\lambda+1)}
	{
	\VStr [\big] \mu
		{\VInter \pi \mu, \VInter X \mu, \VInter Y \mu}
	\ElemSubstr
	\VStr [\big] \kappa {\pi, X, Y}
	},
 \]
where $\pi$ is the function from \Definition {G}, $X = \Set
{\Tuple{p_u, \Beg {p_u}, \End {p_u}}} {u \in \U}$, and $Y = \Set
{\Tuple {v, p_v, \Beg {p_v}, \End {p_v}}} {v \in \U}$.
\end{DEFINITION}

Note that for all $\mu \in \D$, $\beth_\mu = \mu$.

\begin{DEFINITION}{Proj}
For all $\mu \in \D \Union \Braces \kappa$ and $\eta,\nu \in
\Functions \mu 2$ define:
 \ARRAY[lll]{
	\Proj \FComps {\eta,\nu} &=&
	\Set [\Big]  {\FSeq f \in \FComps}
	{\Ind{f_i} \in \ULevel {< \mu} \ForAll i < \SLh f, \\
	&& \Beg f \Subset \eta, \And \End f \Subset \nu}; \\
	\Proj \F {\eta,\nu} &=& \Set [\Big]  {f}
		{\FSeq f \in \Proj \FComps {\eta,\nu}
		\And \SLh f = 1}.
 }
\end{DEFINITION}

\begin{LEMMA}{back_and_forth}
Suppose $\mu \in \D \Union \Braces \kappa$.
\begin{ITEMS}

\ITEM{ind}%
For every $u \in \U$, if $u$ is a limit point or a successor in \UG,
then $p_u \in \Vh \mu$ implies $u \in \ULevel {< \mu}$. For all
successors $u \in \UH$, if $p_u$ is in \Vh \mu then there is $v \in
\ULevel {< \mu}$ such that $p_v = p_u$, $\Beg {p_v} = \Beg {p_u}$, and
$\End {p_v} =
\End {p_u}$.

\ITEM{G}%
For every $v \in \ULevel {< \mu}$ and $\gamma < \mu$, there is $\alpha
\in (\GInds \Inter \mu) \Minus (\gamma +1)$ such that for every
$\eta', \nu' \in \Functions {\beth_\alpha} 2$ with $\Fst v \Subset
\eta'$ and $\Snd v \Subset \nu'$, we can find $u^{\eta',\nu'} \in
\UGLevel \alpha$ satisfying that $\Fst {u^{\eta',\nu'}} = \eta'$,
$\Snd {u^{\eta',\nu'}} = \nu'$, and $u^{\eta',\nu'}$ is a successor of
$v$.

\ITEM{step0}%
If $\eta, \nu \in \Functions \mu 2$ and $q \in \Proj \F {\eta, \nu}$
is such that for $v = \Ind q$ both $\Fst v \Subset \eta$ and $\Snd v
\Subset \nu$ hold, then there is $u \in \ULevel {< \mu}$ satisfying
that $\Ind q \TBelowEq u$ \Note {implying $q \Subset p_u$}, $p_u \in
\Proj \F {\eta, \nu}$, and $\theta \Subset \Dom {p_u} \Inter \Ran
{p_u}$.

\ITEM{step1}%
Suppose $\eta, \nu \in \Functions \mu 2$, $q \in \Proj \F {\eta,
\nu}$, and $\theta < \mu$. There is $\FSeq f \in \Proj \FComps
{\eta,\nu}$ such that $q \Subset f$ and $\theta \Subset \Dom f \Inter
\Ran f$.

\ITEM{step}%
Suppose $\eta, \nu \in \Functions \mu 2$, $\FSeq f \in \Proj \FComps
{\eta,\nu}$, and $\theta < \mu$. There is $\FSeq g \in \Proj \FComps
{\eta,\nu}$ with $g \Superset f$ and $\theta \Subset \Dom g \Inter
\Ran g$.

\end{ITEMS}

\begin{Proof}%
The properties \Item {ind}--\Item {step0} are straightforward
consequences of the definition of the functions $p_u$.
\SimpleComment{%

\ProofOfItem{ind}%
If $u \in \U$ is a limit, then $\beth_{\Ord u} = \Dom {p_u} < \mu$ and
hence $u$ is in \ULevel {< \mu}. If $u$ is a successor in \UG, then
for $\beta = \Ref {\EFunc u}$, $\Dom {p_u} \Union \Ran {p_u} \Subset
\beth_{\beta+1} < \beth_\mu = \mu$. Since $\Ord u = \beta+1$, $u \in
\ULevel {< \mu}$.

Suppose $u$ is a successor in \UH and $p_u \in \Vh \mu$. Since $\Ref
{\EFunc u} < \mu = \beth_\mu$, both \Beg {p_u} and
\End {p_u} are in \Vh \mu. By the definitions of $X$ and $Y$ in
\Definition {D}, for the tuple $\Tuple {p_u, \Beg {p_u}, \End {p_u}}
\in X \Inter \Vh \mu$ there is $v \in \U \Inter \Vh \mu$ such that
\Tuple {v, p_u, \Beg {p_u}, \End {p_u}} is in $Y \Inter \Vh \mu$. So
$v$ is as wanted.

\ProofOfItem{G}%
Suppose first that $v$ is in \Proj \T \TLo \Note {see \Definition
{TLo}}. By \Definition {G}, $\pi$ is a surjective function and the set
\Set {\alpha \in \GInds} {\pi(\alpha) = v} is unbounded in
$\kappa$. By \Definition {D}, also in \Vh \mu, there is $\alpha \in
(\GInds \Inter \mu) \Minus (\gamma +1)$ with $\pi(\alpha) = v$.  Fix
any $\eta'$ and $\nu'$ in \Functions {\beth_\alpha} 2 satisfying $\Fst
v \Subset \eta'$ and $\Snd v \Subset \nu'$. Define $u$ to be the tuple
\Tuple {\eta', \nu', \tau, C} where $C$ is the set $\TCub v \Union
\Braces [\big] {\Ord v}$ and $\tau \in \Functions {\beth_\alpha} 2$ is
such that $\Res \tau \beth_{\Ord v} = \Trd v$ and for all $\xi \in
\beth_\alpha \Minus (\beth_{\Ord v} +1)$, $\tau(\xi) = 0$. Then $u$ is
an element from \UG as wanted. If $v \in \ULevel {< \mu} \Minus \Proj
\T \TLo$, then apply the property proved to \Inv v, and get, for fixed
$\eta', \nu'$ in \Functions {\beth_\alpha} 2, an element
$u^{\eta',\nu'} \in \ULevel {< \mu}$ satisfying the claim for \Inv
v. Then \Inv {(u^{\eta',\nu'})} satisfies the claim for $v$.

\ProofOfItem{step0}%
Each $q' \in \Proj \F {\eta, \nu}$ with $\Ind {q'} = u$ is a partial
function from \Ord u into \Ord u and $\Ord u < \mu$. Hence $\Proj \F
{\eta, \nu} \Subset \Vh \mu$. By \Item {G}, there are $\alpha^1 < \mu$
and $u^1 \in \ULevel {\alpha^1}$ such that $\alpha^1 \geq \Max {\Ord
v, \theta}$, $\Ind q \TBelow u^1$, $\Fst {u^1} = \Res \eta
{\beth_{\alpha^1}}$, and $\Snd {u^1} = \Res \nu {\beth_{\alpha^1}}$.
It might be that $\theta \not\Subset \Dom {p_{u^1}}$ or $\theta
\not\Subset \Ran {p_{u^1}}$. However, by \Item {G} again, there are
$\alpha^2 \in \mu \Minus (\alpha^1+1)$ and $u^2 \in \ULevel
{\alpha^2}$ such that $u^1 \TBelow u^2$, $\Beg {u^2} = \Res \eta
{\beth_{\alpha^2}}$ and $\End {u^2} = \Res \nu {\beth_{\alpha^2}}$.
Then $\Ind q \TBelow u^2$ and $\theta \Subset \Dom {p_{u^2}} \Inter
\Ran {p_{u^2}}$.

}
%
We sketch proofs of the rest of the properties.

\ProofOfItem{step1}%
Here we need the small detail that we used \Id u in \Definition {p}.
Denote \Ind q by $v$. If both $\Fst v \Subset \eta$ and $\Snd v
\Subset \nu$ hold, then the claim follows from \Item {step0}.

Let $\eta', \nu' \in \Functions \mu 2$ be such that $\Fst v \Subset
\eta'$ and $\Snd v \Subset \nu'$. Fix elements $u^0, u^1, u^2$ from
\UG so that
\begin{itemize}

\item $u^0$ is the smallest in $\TBelow$-order with 
$\Fst {u^0} \Subset \eta$, $\Snd {u^0} \Subset \eta'$, and
$\theta \Subset \Dom {p_{u^0}}$;

\item $u^1$ is the smallest in $\TBelow$-order with $\Fst {u^1}
\Subset \eta'$, $\Snd {u^1} \Subset \nu'$, $v \TBelow u^1$ and $\Ran
{\Res {p_{u^0}} \theta} \Subset \Dom {p_{u^1}}$;

\item $u^2$ is the smallest in $\TBelow$-order with $\Fst {u^2}
\Subset \nu'$, $\Snd {u^2} \Subset \nu$, \\
$\Ran {\Res {p_{u^1}} {\Ran {\Res {p_{u^0}} \theta}}} \Subset \Dom
{p_{u^2}}$;

\end{itemize}

Define \FSeq f to be \FuncSeq w W, where $\bar w = \Seq {u^i} {0 \leq
i \leq 2}$. Then \FSeq f is in \Proj \FComps {\eta,\nu}. 

Define $\xi_1$ to be \Min [\zeta \in \Dom q \And \eta(\zeta) \not=
\Fst v (\zeta)] {\zeta+1}, and $\xi_2$ to be \Min [\zeta \in \Dom q
\And \nu(\zeta) \not= \Snd v (\zeta)] {\zeta+1}. Since $\Beg q \Subset
\eta$ and $\End q \subset \nu$, we have that $\xi_1 \geq \PSup {\Dom
q}$ and $\xi_2 \geq \PSup {\Ran q}$. So \EqualRes \eta {\eta'} {\xi_1}
and \EqualRes {\nu'} \nu {\xi_2} ensure that \Res {\Inv {f_0}} {\Dom
q} is identity and \Res {f_2} {\Ran q} is identity. Therefore $q
\Subset f$.

\ProofOfItem{step}%
Since $\Dom f \Union \Ran f$ is bounded in $\mu$ it follows from
\ItemOfLemma {FComps} {limit} that $\Dom {f_i} \Union \Ran {f_i}$ are
bounded in $\mu$ for all $i < \SLh f$. Hence for every $i < \SLh f$,
$f_i \in \Vh \mu$, and by \Item {ind} we may assume, $\Ind {f_i} \in
\Vh \mu$. The claim follows from \Item {step1} by induction on $i <
\SLh f$.
\end{Proof}

\end{LEMMA}
%
\end{SECTION}

\begin{SECTION} {-} {Models} {%
The strongly equivalent non-isomorphic models}
%

Recall that $\kappa$ is a fixed strongly inaccessible cardinal and
$\lambda$ is a fixed regular cardinal below $\kappa$.

For ordinals $\theta < \mu$ and subsets $A$ of $\mu$, \Seqs A \theta
is the set of all $\theta$-sequences of ordinals in $A$. For every
$\theta < \mu < \kappa$ write \BSeqs [\theta] \mu for the set
 \[
	\Set {\BSeq a \in \Seqs \mu \theta}
	{\PSup {\BSeq a} < \mu
	\Text{and for all} i < j < \theta,
	\BSeq [i] a \not= \BSeq [j] a}
 \]
and denote \BigUnion [\theta < \mu] {\BSeqs [\theta] \mu} by \BSeqs
\mu. We write $\Zero$ for the function having domain $\kappa$ and
range \Braces 0.

\begin{DEFINITION}{Models}
For every $\mu \in \D \Union \Braces \kappa$, $\theta < \mu$, and
$\BSeq a \in \BSeqs [\theta] \mu$ we define a family \Seq {\Rel [\eta]
a} {\eta \in \Functions \mu 2} of relations \Note {having arity
$\theta$} on $\mu$ as follows:

relations \Rel [\eta] a, $\eta \in \Functions \mu 2$, are the smallest
subsets of \Seqs \mu \theta having the properties:
\begin{itemize}

\item if $\eta = \Zero [\mu]$ then $\BSeq a \in \Rel [\eta] a$;

\item if there is $\FSeq f \in \Proj \FComps {\Zero [\mu], \eta}$ with
$\Dom f = \BSeq a$, then $f(\BSeq a) \in \Rel [\eta] a$.

\end{itemize}
Suppose $\mu \in \D \Union \Braces \kappa$. Define \Voc \mu to be the
vocabulary \Set {\Rel a} {\BSeq a \in \BSeqs \mu} where each \Rel a is
a relation symbol of arity \BLh a. For every $\eta \in \Functions \mu
2$, let \Model \eta be the \Voc \mu-structure with domain $\mu$ and
interpretations $\Int {(\Rel a)} {\Model \eta} = \Rel [\eta] a$ for
all $\BSeq a \in \BSeqs \mu$. For all $\chi \in \D \Inter \mu$ and $A
\Subset \mu$, we write \ModelRes \eta A \chi for the models having
vocabulary \Voc \chi, domain $A$, and interpretations $\Int {(\Rel a)}
{\ModelRes \eta A \chi} = \Rel [\eta] a \Inter \Seqs A {\BLh a}$ for
all $\BSeq a \in \BSeqs \chi$.
\end{DEFINITION}

\begin{FACT}{Relations}%
Assume $\mu \in \D \Union \Braces \kappa$ and $\eta \in \Functions \mu
2$.
\begin{ITEMS}

\ITEM{bounded}%
For every $\BSeq a \in \BSeqs \mu$, \Rel [\eta] a is a subset of
\BSeqs \mu.

\ITEM{res}%
For all $\chi \in \D \Inter \mu$, $\Model {\Res \eta \chi} = \ModelRes
\eta \chi \chi$.

\end{ITEMS}

\begin{Proof}%
\ProofOfItem{bounded}%
Assume that for some $\BSeq b \in \Rel [\eta] a$, $\PSup {\BSeq b} =
\mu$. Then there should be $\FSeq f \in \Proj \FComps {\Zero
[\mu], \eta}$ with $\Dom f = \BSeq a$ and $f(\BSeq a) = \BSeq b$
contrary to \ItemOfLemma {FComps} {bounded} and the fact $\PSup {\BSeq
a} < \mu$. %
\PvComment {This property also follows directly from \ItemOfLemma
{back_and_forth} {ind}.}

\ProofOfItem{res}%
Abbreviate \Res \eta \chi by $\nu$ and let \BSeq a be a sequence from
\BSeqs \chi. The interpretation $\Int {(\Rel a)} {\Model \nu} = \Rel
[\nu] a$ is a subset of the interpretation \Int {(\Rel a)} {\ModelRes
\eta \chi \chi} since $\Proj \FComps {\Zero [\chi], \nu} \Subset \Proj
\FComps {\Zero [\mu], \eta}$. Suppose $\BSeq b \in \Int {(\Rel a)}
{\ModelRes \eta \chi \chi}$ and let $\FSeq f \in \Proj \FComps {\Zero
[\mu], \eta}$ be such that $\Dom f = \BSeq a$ and $f(\BSeq a) = \BSeq
b$. By \ItemOfLemma {back_and_forth} {ind}, we may assume $\Ind {f_i}
\in \ULevel {<\chi}$ for every $i < \SLh f$. Consequently, $\FSeq f
\in \Proj \FComps {\Zero [\mu'], \nu}$ and $\BSeq b \in \Rel [\nu] a$.
\end{Proof}%
\end{FACT}

\begin{FACT}{Aux_Isom}%
Suppose $\mu \in \D \Union \Braces \kappa$ and $\eta, \nu \in \Functions
\mu 2$.
\begin{ITEMS}

\ITEM{partial}%
For all $v \in \U$ with $\Fst v \Subset \eta$ and $\Snd v \Subset
\nu$, the function $p_v$ is a partial isomorphism from \Model \eta
into \Model \nu.

\ITEM{rel}
For every $\theta < \mu$ and $\BSeq b \not= \BSeq c \in \BSeqs
[\theta] \mu$, if there exist $\BSeq a \in \BSeqs [\theta] \mu$
satisfying
\[
		\Model \eta \models \Rel a (\BSeq b) \And
		\Model \nu \models \Rel a (\BSeq c),
\]
then there is $\FSeq f \in \Proj \FComps {\eta,\nu}$ with
$f(\BSeq b) = \BSeq c$.

\end{ITEMS}

\begin{Proof}%
Both of these properties are direct consequences of \Definition {Models} and
\Fact {Minimal}.
\end{Proof}%
\end{FACT}

\begin{LEMMA}{Characterization}
For all $s, t \in \Functions \kappa 2$, 
\begin{ITEMS}

\ITEM{equiv->isom}
$s \Equiv t$ implies $\Model s \Isomorphic \Model t$ \Note {$\Equiv$
is given in \Definition {Sigma}}, and

\ITEM{isom->equiv}
if $\Model s \Isomorphic \Model t$ then $s \Equiv t$.

\end{ITEMS}

%
\begin{Proof}%
\ProofOfItem{equiv->isom}
Suppose $s \Equiv t$, and let \Function r \kappa 2 be such that $\VStr
\kappa {P, s, t, r} \models \phi$. For every $\delta \in C' = \Cub str
\Inter \D$ define $u_\delta$ to be the tuple \Tuple {\Res s \delta,
\Res t \delta, \Res r \delta, \Cub str \Inter \delta} \Note {\Cub str
is given in \Definition {H} and \D is given in \Definition
{D}}. Directly by \Definition {H}, for all $\delta < \epsilon \in C'$,
$u_\delta, u_\epsilon$ are in \UH and $u_\delta \TBelow u_\epsilon$.
Hence $p_{u_\delta} \Subset p_{u_\epsilon}$ for $\delta < \epsilon \in
C'$, and moreover, for the function $h = \BigUnion [\delta \in C']
{p_{u_\delta}}$ both of the equations $\Dom h = \kappa$ and $\Ran h =
\kappa$ hold.  Consequently $h$ is an isomorphism from \Model s onto
\Model t.

\ProofOfItem{isom->equiv}
Suppose $s \not= t$ and for fixed $\xi < \kappa$, $s(\xi) \not=
t(\xi)$.  Let $h$ be an isomorphism from \Model s onto \Model t, and
let $S'$ be the set \Set {\delta \in \kappa \Minus (\xi+1)} {h[\delta]
= \delta \Text {is a cardinal of cofinality}\, \geq \lambda}.  Since
$h$ is an isomorphism and $\Res s \delta \not= \Res t \delta$ for all
$\delta \in S'$, it follows from \ItemOfFact {Aux_Isom} {rel} that for
every $\delta \in S'$ there is a sequence ${\FSeq f}^\delta \in \Proj
\FComps {s, t}$ such that $f^\delta = \Res h \delta$. For all $\delta
< \epsilon \in S'$, $f^\delta = \Res h \delta \Subset \Res h \epsilon
= f^\epsilon$. Since $S'$ is stationary in $\kappa$, there are $n <
\omega$ and a stationary subset $S$ of $S'$ such that the equation
$\SLh [^\delta] f = n$ holds for every $\delta \in S$.

Consider some $\delta \in S$ and $i < n$. Abbreviate \Ind {f^\delta_i}
by $u^\delta_i$. By \ItemOfLemma {FComps} {bounded}, $\Dom
{f^\delta_i} = \Ran {f^\delta_i} = \delta$.  Moreover, by \ItemOfLemma
{Limit} {componentwise}, $u^\delta_i \TBelowEq u^\epsilon_i$ and
$f^\delta_i \Subset f^\epsilon_i$ for all $\epsilon \in S \Minus
\delta$. By \ItemOfFact {p} {cofinality}, $u^\delta_i$ is in \UH.  So
$f^\delta_i = p_{u^\delta_i}$ and $\Dom {f^\delta_i} = \delta = \Ord
{u^\delta_i} = \beth_\delta$.  Define for each $i < n$,
 \ARRAY{
	s_i = \BigUnion [\delta \in S] {
		\Fst {u^\delta_i}
	}; \\
	r_i = \BigUnion [\delta \in S] {
		\Trd {u^\delta_i}
	};
 }
and let $s_n$ be \BigUnion [\delta \in S] {\Snd
{u^\delta_{n-1}}}. Then $s = s_0$ and $t = s_n$.

We claim that $s \Equiv t$. By the transitivity of $\Equiv$ it is
enough to show that for every $i < n$, $r_i$ witness $s_i \Equiv
s_{i+1}$. Contrary to this subclaim assume that for some $i < n$,
 \[
	\VStr \kappa {P, s_i, s_{i+1}, r_i} \not\models \phi.
 \]
Then there is $\delta \in S$ for which
 \[
	\VStr \delta {\VInter P \delta,
	\TripleRes {s_i} {s_{i+1}} {r_i} \delta}
	\ElemSubstr
	\VStr \kappa {P, s_i, s_{i+1}, r_i}.
 \]
However $\Res {s_i} \delta = \Fst {u^\delta_i}$, $\Res {s_{i+1}}
\delta = \Snd {u^\delta_i}$, and $\Res {r_i} \delta = \Trd
{u^\delta_i}$, and so
 \[
	\VStr \delta {\VInter P \delta,
	\Fst {u^\delta_i}, \Snd {u^\delta_i}, \Trd {u^\delta_i}}
	\not\models \phi,
 \]
contrary to the fact that $u^\delta_i$ is in \UH.
\end{Proof}%
%

\end{LEMMA}

In the following two lemmas we assume existence of a regular cardinal
$\mu$ in \D. Such $\mu$ does not necessarily exists, if $\kappa$ is an
arbitrary strongly inaccessible cardinal. However, these lemmas are
only preliminaries for the main lemma, \Lemma {No(M)<2^k}, where we
assume $\kappa$ to be a weakly compact cardinal. Note, when $\mu =
\kappa$ in \ItemOfLemma {Models} {equivalent} below, it suffices that
$\kappa$ is a strongly inaccessible cardinal.

\begin{LEMMA}{Models}%
Suppose $\mu \in \D$ is a regular cardinal or $\mu = \kappa$, and that
$\eta, \nu$ are functions from $\mu$ into $2$.
\begin{ITEMS}

\ITEM{equivalent}
$\Model \eta \MEquiv \mu \lambda \Model \nu$.

\ITEM{in_some_rel}%
For every $\theta < \mu$, the model \Model \eta satisfies the \Lan
\mu-sentence
 \[
	\forall \Seq {x_i} {i < \theta} \Par [\Big] {
		\bigvee_{\BSeq a \in \BSeqs [\theta] \mu}
		\Rel a \Par[\big] {\Seq {x_i} {i < \theta}}
	}
 \]

\ITEM{ext}%
For all $\BSeq a \in \BSeqs \mu$ and $\xi < \mu$, the following \Lan
\mu-sentence holds in \Model \eta:
 \[
	\forall \bar x \Par [\Big] {
		\Rel a (\bar x) \rightarrow 
		\exists y \Par [\big] {
			\Relation {\Concat {\SimpleSeq \xi} {\BSeq a}}
			(\Concat {\SimpleSeq y} {\bar x})
		}
	}.
 \]

\ITEM{ext_rel}%
For all $\BSeq a \in \BSeqs \mu$, \Model \eta satisfies the \Lan
\mu-sentence:
 \[
	\forall \bar x \forall y \Par [\Big] {
		\Rel a (\bar x) \rightarrow
		\bigvee_{\xi < \mu}
		\Relation {\Concat {\SimpleSeq \xi} {\BSeq a}}
		(\Concat {\SimpleSeq y} {\bar x})
	}.
 \]

\end{ITEMS}

\begin{Proof}%
\ProofOfItem{equivalent}%
We give a winning strategy for player \PlayerTwo in the game \EF \mu
\lambda {\Model \eta} {\Model \nu} \Note {see \Definition
{EF}}. Suppose $i < \lambda$ and for each $j \leq i$, player
\PlayerOne has chosen $X_j \in \Braces {\Model \eta, \Model \nu}$ and
$A_j \Subset \mu$ \Note {where $\mu$ is the domain of both \Model \eta
and \Model \nu}. Suppose that for every $j < i$, player \PlayerTwo has
replied with a partial isomorphism $p_{u_j}$ satisfying that $u^j \in
\UG$, $\Fst {u^j} \Subset \eta$, $\Snd {u^j} \Subset \nu$, $\BigUnion
[k \leq j] {A_k} \Subset \Dom {p_{u_j}} \Inter \Ran {p_{u_j}}$, and
for all $k < j$, $u^k \TBelow u^j$. Since $i < \lambda$ and $u^j \in
\UG$ for each $j < i$, the tuple $v = \BigUnion [j < i] {u^j}$ is in
\UG by \ItemOfFact {U} {G-closed}. Let $\theta$ be the smallest
ordinal which is strictly greater than any ordinal in \BigUnion [j
\leq i] {A_j} \Note {$\theta < \mu$ since $\mu$ is regular, $i < \mu$,
and $\Card {A_j} < \mu$ for every $j \leq i$}. By \ItemOfLemma
{back_and_forth} {step1}, there is $u^i$ in \ULevel {< \mu} satisfying
that $v \TBelow u^i$, $\Fst {u^i} \Subset \eta$, $\Snd {u_i} \Subset
\nu$, and $\theta \Subset \Dom {p_{u^i}} \Inter \Ran {p_{u^i}}$. Since
$\BigUnion [j < i] {p_{u^j}} = p_v \Subset p_{u^i}$, the partial
isomorphism $p_{u^i}$ is a valid reply for player \PlayerTwo in the
round $i$.

\ProofOfItem{in_some_rel}%
By \Definition {Models}, for every $b \in \BSeqs [\theta] \mu$, $\Rel
b(\BSeq b)$ is satisfied in \Model {\Zero [\mu]}. The claim follows
from \Item {equivalent}.

\ProofOfItem{ext}%
By \Item {equivalent} we may assume $\eta = \Zero [\mu]$. For $\bar x
= \BSeq a$ the claim holds directly by \Definition {Models}. For any
$\bar x = \BSeq b \in \Rel [\Zero] a \Minus \Braces {\BSeq a}$ there
is some $\FSeq f \in \Proj \FComps {\Zero [\mu], \Zero [\mu]}$ such
that $\Dom f = \BSeq a$ and $f(\BSeq a) = \BSeq b$. Since $\mu \in \D$
there is by \ItemOfLemma {back_and_forth} {step}, $\FSeq g \in \Proj
\FComps {\Zero [\mu], \Zero [\mu]}$ with $g \Superset f$ and $\Dom g =
\Concat {\SimpleSeq \xi} {\BSeq a}$.

\ProofOfItem{ext_rel}%
Analogously to the proof of \Item {ext}. If $\bar x = \BSeq b$ and $y
= \zeta$ then there is some $\FSeq f \in \Proj \FComps {\Zero
[\mu], \Zero [\mu]}$ such that $f(\BSeq a) = \BSeq b$. Moreover by
\ItemOfLemma {back_and_forth} {step}, there is $\FSeq g \in \Proj
\FComps {\Zero [\mu], \Zero [\mu]}$ with $g \Superset f$ and $\Ran g =
\Concat {\SimpleSeq \zeta} {\BSeq b}$.
\end{Proof}
\end{LEMMA}

\begin{LEMMA}{IsomorphismType}%
Suppose $\mu$ is a regular cardinal in \D and $A$ is a subset of
$\kappa$ having cardinality $\mu$.%
\begin{ITEMS}%

\ITEM{aux}%
Suppose $\eta \in \Functions \mu 2$, $A \Subset \mu$, and $\ModelRes
\eta A \mu \LEquiv \mu \Model \eta$. Then $A = \mu$.

\ITEM{main}%
If $\ModelRes \Zero A \mu \LEquiv \mu \Model {\Zero [\mu]}$ then there
is $\eta \in \Functions \mu 2$ for which $\ModelRes \Zero A \mu
\Isomorphic \Model \eta$.

\end{ITEMS}

%
\begin{Proof}%
\ProofOfItem{aux}%
Suppose, contrary to the claim, that $\xi < \mu$ is not in $A$. Let
\BSeq b in \BSeqs [\omega] \mu be such that $\BSeq b \Subset A$,
$\BSeq [0] b > \xi$ and for every $i < \omega$, $\BSeq [i] b < \Card
{\BSeq [i+1] b}$. This is possible since $\mu = \Cf \mu$, $\mu \in \D$
implies $\mu$ is an uncountable limit cardinal, and $\Card A = \mu$
implies that $A$ is unbounded in $\mu$.  By \ItemOfLemma {Models}
{in_some_rel}, there is $\BSeq a \in \BSeqs [\omega] \mu$ such that
$\Model \eta \models \Rel a (\BSeq b)$. Since $\BSeq b \Subset A$ also
\ModelRes \eta A \mu satisfies $\Rel a(\BSeq b)$. By \ItemOfLemma
{Models} {ext_rel}, there is $\xi' < \mu$ such that $\Model \eta
\models \Relation {\Concat {\SimpleSeq {\xi'}} {\BSeq a}} (\Concat
{\SimpleSeq \xi} {\BSeq b})$. By \ItemOfLemma {Models} {ext}, there
should be $\zeta \in A$ with $\ModelRes \eta A \mu \models \Relation
{\Concat {\SimpleSeq {\xi'}} {\BSeq a}} (\Concat {\SimpleSeq {\zeta}}
{\BSeq b})$. However, then by \ItemOfFact {Aux_Isom} {rel}, there
should be $\FSeq f \in \FComps$ satisfying $\Dom f = \Braces \xi
\Union \BSeq b$, $f(\xi) = \zeta \not= \xi$, and $f(\BSeq b) = \BSeq
b$ contrary to \ItemOfLemma {Limit} {id}.

\ProofOfItem{main}%
Our proof has the following structure:
\begin{itemize}

\item %
When $A \Subset \mu$ the claim follows from \Item {aux} and
\ItemOfFact {Relations} {res}.

\item %
The case that $A$ is not a subset of $\zeta + \mu$ for any $\zeta \in
A$ is shown to be impossible.

\item %
Lastly we prove that when $A \Subset \zeta + \mu$ for some $\zeta \in
A$, there are $\eta \in \Functions \mu 2$ and $\FSeq g \in \FComps$
such that $\Dom g = \mu$, $\Ran g = A$, $\Beg g = \eta$, and $\End g
\Subset \Zero$. So $g$ is an isomorphism between \Model \eta and
\ModelRes \Zero A \mu.

\end{itemize}

Suppose there is an $\omega$-sequence $\BSeq b$ such that $\BSeq [0] b
> \mu$ and for all $l < \omega$, $\BSeq [l] b \in A$ and $\BSeq [l+1]
b \geq \BSeq [l] b + \mu$. By the equivalence $\ModelRes \Zero A \mu
\LEquiv \mu \Model {\Zero [\mu]}$ and \ItemOfLemma {Models}
{in_some_rel}, there is $\BSeq a \in \BSeqs [\omega] \mu$ such that
$\ModelRes \Zero A \mu \models \Rel a (\BSeq b)$. Hence there should
be $\FSeq f \in \FComps$ with $\Dom f = \BSeq a$ and $f(\BSeq a) =
\BSeq b$ contrary to \ItemOfLemma {Ran} {limit}.

Suppose $\zeta \in A$ and $A \Subset \zeta + \mu$. As above, there are
$\FSeq f \in \FComps$ and $\gamma < \mu$ with $\Dom f = \Braces
\gamma$, $f(\gamma) = \zeta$, $\Beg f = \Res \Zero {(\gamma +1)}$, and
$\End f = \Res \Zero {(\zeta +1)}$.  Since $\gamma < \mu \leq \zeta =
f(\gamma)$, there is the smallest index $k < \SLh f$ such that
$\ResComp f {< k} (\gamma) < \mu$ and $\ResComp f {\leq k} (\gamma)
\geq \mu$. By \ItemOfLemma {FComps} {increasing}, $f_j$ is increasing
for all $j \in \Interval k {\SLh f-1}$. Let $\bar u$ be the sequence
\Seq {\Ind {f_j}} {j \in \Interval k {\SLh f -1}}. By \ItemOfLemma
{Ran} {seq}, the sequence \FuncSeq u \mu is a well-defined member of
\FComps, and moreover, $\End {\Func u \mu} \Subset \Zero$. Abbreviate
this sequence by $g$ and the ordinal $\ResComp f {< k} (\gamma)$ by
$\xi$. We define the wanted $\eta$ to be \Beg g.

Finally we show that for this $g$ the equation $A = \Ran g$
holds. Suppose $\zeta'$ is in $A$ but not in \Ran g. By the
equivalence and \ItemOfLemma {Models} {in_some_rel}, there are
$\epsilon, \epsilon' < \mu$ such that $R_{\SimpleSeq \epsilon}(\zeta)$
and $R_{\SimpleSeq {\epsilon', \epsilon}}(\zeta', \zeta)$ hold in
\ModelRes \Zero A \mu. By \ItemOfLemma {Models} {ext}, there is $\xi'
< \mu$ for which $R_{\SimpleSeq {\epsilon', \epsilon}} (\xi', \xi)$ in
\Model \eta. However, by \ItemOfFact {Aux_Isom} {rel}, there should be
$\FSeq h \in \FComps$ with $h(\xi) = \zeta$ and $h(\xi') = \zeta'$,
contrary to \ItemOfLemma {Ran} {ran}. From another direction, if $A
\ProperSubset \Ran g$, then $\eta$ and the set $B = \Image {\Inv g} A
\ProperSubset \mu$ contradict \Item {aux}, since by \ItemOfLemma
{Models} {equivalent}, $\Model \eta \LEquiv \mu \Model {\Zero [\mu]}$,
by our assumption, $\Model {\Zero [\mu]} \LEquiv \mu \ModelRes \Zero A
\mu$, and \Isomorphism {\Res {\Inv g} A} {\ModelRes \Zero A \mu}
{\ModelRes \eta B \mu}.
\end{Proof}%
%

\end{LEMMA}

\begin{LEMMA}{No(M)<2^k}%
Suppose $\kappa$ is a weakly compact cardinal and \M is a model of
cardinality $\kappa$ with $\M \LEquiv \kappa \Model \Zero$. Then there
is $s \in \Functions \kappa 2$ for which $\M \Isomorphic \Model s$.

\begin{Proof}%
Without loss of generality we may assume that the domain of \M is
$\kappa$. By the $\LEquiv \kappa$-equivalence of the models \M and
\Model \Zero, let for every regular cardinal $\mu < \kappa$, $A_\mu$
be a subset of $\kappa$ such that $\Res {\VocRes \M \mu} \mu
\Isomorphic \ModelRes \Zero {A_\mu} \mu$. Let $Y$ be the set given in
\Definition {D}. Note that for all $\mu \in \D \Union \Braces \kappa$
and $\eta \in \Functions \mu 2$, the model \Model \eta is definable
from $\eta$ and $Y \Inter \Vh \mu$. Let $\tau$ be a winning strategy
for \PlayerTwo in the game \EF \kappa \omega \M {\Model \Zero}. Assume
now, contrary to the claim, that $\M \not\Isomorphic \Model s$ for all
$s \in \Functions \kappa 2$. Because $\kappa$ is
$\Pi^1_1$-indescribable, there is a regular cardinal $\mu < \kappa$
such that \VStr \mu {\Res {\VocRes \M \mu} \mu, \VInter \tau \mu,
\VInter Y \mu} satisfies the following:
 \[
	\ForAll \eta \in \Functions \mu 2,
	\ModelRes {} \mu \mu \not\Isomorphic \Model \eta.
 \]
Then $\Model {\Zero [\mu]} = \ModelRes \Zero \mu \mu \LEquiv \mu \Res
{\VocRes \M \mu} \mu$, and by the isomorphism $\Res {\VocRes \M \mu}
\mu \Isomorphic \ModelRes \Zero {A_\mu} \mu$, we have $\Model {\Zero
[\mu]} \LEquiv \mu \ModelRes \Zero {A_\mu} \mu$ and for all $\eta \in
\Functions \mu 2$, $\ModelRes \Zero {A_\mu} \mu \not\Isomorphic \Model
\eta$.  This contradicts \ItemOfLemma {IsomorphismType} {main}.
\end{Proof}%
\end{LEMMA}

\begin{LEMMA}{No}%
Suppose $\kappa$ is a weakly compact cardinal, $\lambda < \kappa$ is a
regular cardinal, and there is a $\Sigma^1_1$-equivalence relation on
\Functions \kappa 2 having $\mu$ different equivalent classes. Then
there exists a model \M such that the vocabulary of \M consists of one
relation symbol of finite arity, $\Card \M = \kappa$, and $\No
[\lambda] \M = \mu$.

\begin{Proof}%
By the preceding lemmas the model \Model \Zero defined as in
\Definition {Models} satisfies the claim, except that the vocabulary
of \M is overly large. However, by \cite [Claim 1.3(1)] {Sh189}, the
inaccessibility of $\kappa$ ensures that there is a model \N of
cardinality $\kappa$ with $\lambda$ many relations of finite arity
satisfying $\No N = \No {\Model \Zero}$ \Note {the proof is a simple
coding}. Furthermore, by \cite [Claim 1.4(2)] {Sh189}, the $\lambda$
many relations can be coded by one relation so that the other
properties are preserved.  Actually, the claims cited concern the case
$\lambda = \aleph_0$, but there is no problem to preserve $\No
[\lambda] {\Model \Zero}$ in the cases $\aleph_0 < \lambda < \kappa$,
too.
\end{Proof}
\end{LEMMA}%
\end{SECTION}

%

%

%

\begin{tabbing}
Saharon Shelah:\= \\
	\>Institute of Mathematics\\
	\>The Hebrew University\\
	\>Jerusalem. Israel\\
\\
	\>Rutgers University\\
	\>Hill Ctr-Busch\\
	\>New Brunswick. New Jersey 08903\\
	\>\texttt{shelah@math.rutgers.edu}
\end{tabbing}

\begin{tabbing}
Pauli V\"{a}is\"{a}nen:\= \\
	\>Department of Mathematics\\
	\>P.O. Box 4\\
	\>00014 University of Helsinki\\
	\>Finland\\
	\>\texttt{pauli.vaisanen@helsinki.fi}
\end{tabbing}

%

\end{document}